\documentclass[aop,MSNbibl,seceqn,nameyear,dvips]{arximspdf}
\usepackage{mathrsfs}

\doi{10.1214/12-AOP759} %kopijuoti is PTS
\volume{42}
\issue{1}
\pubyear{2014}
\firstpage{41}
\lastpage{79}

\makeatletter

\newtheorem{theorem}{Theorem}[section]
\newtheorem{prop}[theorem]{Proposition}

\newcommand{\iittTheta}{{\Theta}}
\newcommand{\iittOmega}{{\Omega}}
\newcommand{\iittPsi}{{\Psi}}

\newcommand{\mcr}{\mathscr}
\newcommand{\mbb}{\mathbb}
\newcommand{\mbf}{\mathbf}
\newcommand{\mrm}{\mathrm}

\makeatother

\begin{document}
\begin{frontmatter}

\title{Path-valued branching processes and nonlocal branching superprocesses\thanksref{T1}}
\runtitle{Path-valued processes and superprocesses}

\thankstext{T1}{Supported by NSFC (No. 11131003), 973 Program (No. 2011CB808001) and 985 Program.}

\begin{aug}
\author[A]{\fnms{Zenghu} \snm{Li}\corref{}\ead[label=e1]{lizh@bnu.edu.cn}\ead[label=u1,url]{http://math.bnu.edu.cn/\textasciitilde lizh/}}
\runauthor{Z. Li}
\affiliation{Beijing Normal University}
\address[A]{School of Mathematical Sciences\\
Beijing Normal University\\
Beijing 100875\\
P. R. China\\
\printead{e1}\\
\printead{u1}} %adresu isvedimo komanda gale!
\end{aug}

% HISTORY:
\received{\smonth{11} \syear{2011}}
\revised{\smonth{3} \syear{2012}}

% ABSTRACT
%
\begin{abstract}
A family of continuous-state branching
processes with immigration are constructed as the solution flow of a
stochastic equation system driven by time--space noises. The family can be
regarded as an inhomogeneous increasing path-valued branching process
with immigration. Two nonlocal branching immigration superprocesses can
be defined from the flow. We identify explicitly the branching and
immigration mechanisms of those processes. The results provide new
perspectives into the tree-valued Markov processes of Aldous and Pitman
[\textit{Ann. Inst. Henri Poincar\'e Probab. Stat.} \textbf{34} (1998)
637--686] and Abraham and Delmas
[\textit{Ann. Probab.} \textbf{40} (2012) 1167--1211].
\end{abstract}

% KEYWORDS

\end{frontmatter}

%s1 #&#
\section{Introduction}\label{sec1}

Continuous-state branching processes (CB-processes) are positive Markov
processes introduced by \citet{Jir58} to model the evolution of
large populations of small particles. Continuous-state branching
processes with immigration (CBI-processes) are generalizations of them
describing the situation where immigrants may come from other sources
of particles; see, for example, \citet{KawWat71}. The law of a
CB-process is determined by its \textit{branching mechanism} $\phi$,
which is a function with the representation
%
%e1.1 #&#
\begin{equation}
\label{1.1} \phi(\lambda) = b\lambda+ \frac{1}{2}\sigma^2
\lambda^2 + \int_0^\infty
\bigl(e^{-z\lambda}-1+z\lambda\bigr)m(dz),
\end{equation}
where $\sigma\ge0$ and $b$ are constants, and $(z\land z^2)m(dz)$ is a
finite measure on $(0,\infty)$. In most cases, we only define the
function $\phi$ on $[0,\infty)$, but it can usually be extended to an
analytic function on an interval strictly larger than $[0,\infty)$. The
branching mechanism is said to be \textit{critical}, \textit{subcritical}
or \textit{supercritical} according as $b=0$, $b>0$ or $b<0$.

A CB-process can be obtained as the small particle limit of a sequence
of discrete Galton--Watson branching processes; see, for example,
\citet{Lam67}. A~genealogical tree is naturally associated with a
Galton--Watson process. The genealogical structures of CB-processes
were investigated by introducing \textit{continuum random trees} in the
pioneer work of Aldous (\citeyear{Ald91}, \citeyear{Ald93}),
where the \textit{quadratic branching mechanism} $\phi(\lambda) =
\lambda^2$ was considered. Continuum random trees corresponding to
general branching mechanisms were constructed in Le~Gall and Le~Jan
(\citeyear{LeGLeJ98N1,LeGLeJ98N2}) and were studied further in
\citet{DuqLeG02}. By pruning a Galton--Watson tree,
\citet{AldPit98} constructed a decreasing tree-valued process.
Then they used time-reversal to obtain an increasing tree-valued
process starting with the trivial tree. They gave some
characterizations of the increasing process up to the \textit{ascension
time}, the first time when the increasing tree becomes infinite.

Tree-valued processes associated with general CB-processes were studied
in \citet{AbrDel12}. By shifting a critical branching
mechanism, they defined a family of branching mechanisms
$\{\psi_\theta\dvtx  \theta \in \iittTheta\}$, where $\iittTheta=
[\theta_\infty,\infty)$ or $(\theta_\infty,\infty)$ for some
$\theta_\infty\in[-\infty,0]$. \citet{AbrDel12} constructed a
decreasing tree-valued Markov process $\{\mcr{T}_\theta\dvtx
\theta\in\iittTheta\}$ by pruning a continuum tree, where the tree
$\mcr{T}_\theta$ has branching mechanism $\psi_\theta$. The
\textit{explosion time} $A$ was defined as the smallest negative time
when the tree (or the total mass of the corresponding CB-process) is
finite. \citet{AbrDel12} gave some characterizations of the
evolution of the tree after this time under an excursion law. For the
quadratic branching mechanism, they obtained explicit expressions for
some interesting distributions. Those extend the results of
\citet{AldPit98} on Galton--Watson trees in the time-reversed
form. The main tool of \citet{AbrDel12} was the exploration
process of Le~Gall and Le~Jan (\citeyear{LeGLeJ98N1,LeGLeJ98N2}) and
\citet{DuqLeG02}. Some general ways of pruning random trees
in discrete and continuous settings were introduced in
\citet{AbrDelHe}, \citet{AbrDelVoi10}.

In this paper, we study a class of increasing path-valued Markov
processes using the techniques of stochastic equations and measure-valued
processes developed in recent years. Those path-valued processes are
counterparts of the tree-valued processes of \citet{AbrDel12}. A
special case of the model is described as follows. Let $T = [0,\infty)$
or $[0,a]$ or $[0,a)$ for some $a>0$. Let $(\theta,\lambda)\mapsto
\zeta_\theta(\lambda)$ be a continuous function on $T\times
[0,\infty)$
with the representation
\[
\zeta_\theta(\lambda) = \beta_\theta\lambda+ \int
_0^\infty\bigl(1-e^{-z\lambda}\bigr)
n_\theta(dz),\qquad \theta\in T, \lambda\ge0,
\]
where $\beta_\theta\ge0$ and $zn_\theta(dz)$ is a finite kernel
from $T$ to
$(0,\infty)$. Let $\phi$ be a branching mechanism given by (\ref{1.1}).
Under an integrability condition, the function
%
%e1.2 #&#
\begin{equation}
\label{1.2} \phi_q(\lambda):= \phi(\lambda) - \int
_0^q \zeta_\theta(\lambda)\,d\theta,\qquad
\lambda\ge0,
\end{equation}
also has the representation (\ref{1.1}) with the parameters $(b,m) =
(b_q,m_q)$ depending on $q\in T$. Let $m(dy,dz)$ be the measure on
$T\times(0,\infty)$ defined by
\[
m\bigl([0,q]\times[c,d]\bigr) = m_q[c,d],\qquad q\in T, d>c>0.
\]
Let $W(ds,du)$ be a white noise on $(0,\infty)^2$ based on the Lebesgue
measure, and let $\tilde{N}_0(ds,dy,dz,du)$ be a compensated Poisson random
measure on $(0,\infty)\times T\times(0,\infty)^2$ with intensity
$dsm(dy,dz)\,du$. Let $\mu\ge0$ be a constant. For $q\in T$, we
consider the
stochastic equation
%
%e1.3 #&#
\begin{eqnarray}
\label{1.3} X_t(q) &=& \mu- b_q\int_0^t
X_{s-}(q) \,ds + \sigma\int_0^t\int
_0^{X_{s-}(q)} W(ds,du)
\nonumber\\[-8pt]\\[-8pt]
&&{} + \int_0^t\int_{[0,q]}
\int_0^\infty\int_0^{X_{s-}(q)}
z \tilde{N}_0(ds,dy,dz,du).
\nonumber
\end{eqnarray}
We shall see that there is a pathwise unique positive c\`adl\`ag solution
$\{X_t(q)\dvtx\break   t\ge0\}$ to (\ref{1.3}). Then we can talk about the solution
flow $\{X_t(q)\dvtx  t\ge0, q\in T\}$ of the equation system. We prove that
each $\{X_t(q)\dvtx  t\ge0\}$ is a CB-process with branching mechanism
$\phi_q$, and $\{(X_t(q))_{t\ge0}\dvtx  q\in T\}$ is an inhomogeneous
path-valued increasing Markov process with state space $D^+[0,\infty)$,
the space of positive c\`adl\`ag paths on $[0,\infty)$ endowed with the
Skorokhod topology.

The formulation of path-valued processes provides new perspectives into
the evolution of the random trees of \citet{AldPit98} and
\citet{AbrDel12}. From this formulation we can derive some structural
properties of the model that have not been discovered before. For $q\in
T$ let us define the random measure $Z_q(dt) = X_t(q)\,dt$ on $[0,\infty)$.
We shall see that $\{Z_q\dvtx  q\in T\}$ is an inhomogeneous increasing
superprocess involving a nonlocal branching structure, and the total mass
process
\[
\sigma(q):= \int_0^\infty X_s(q)\,ds,\qquad
q\in T,
\]
is an inhomogeneous CB-process. Then one can think of $\{X(q)\dvtx  q\in T\}$
as a path-valued branching process. On the other hand, for each $t\ge0$
the random increasing function $q\mapsto X_t(q)$ induces a random measure
$Y_t(dq)$ on $T$ such that $X_t(q) = Y_t[0,q]$ for $q\in T$. We prove
that $\{Y_t\dvtx  t\ge0\}$ is a homogeneous superprocess with both local and
nonlocal branching structures. We also establish some properties of an
excursion law $\mbf{N}_0$ for the superprocess $\{Y_t\dvtx  t\ge0\}$.
Given a
branching mechanism $\phi$ of the form (\ref{1.1}), for a suitable
interval $T$ we can define a family of branching mechanisms $\{\phi_q\dvtx
q\in T\}$ by
\[
\phi_q(\lambda) = \phi(\lambda-q) - \phi(-q),\qquad \lambda\ge0,
\]
where the two terms on the right-hand side are defined using (\ref{1.1}).
The family can be represented by (\ref{1.2}) with $\zeta_\theta
(\lambda) =
- (\partial/\partial\lambda) \phi_\theta(\lambda)$. In this case, the
path-valued process $q\mapsto(X_t(q))_{t\ge0}$ under the excursion law
$\mbf{N}_0$ corresponds to the time-reversal of the tree-valued process
$\theta\mapsto\mcr{T}_\theta$ of \citet{AbrDel12}. In general,
we may associate $\{X(q)\dvtx  q\in T\}$ with a ``forest-valued branching
process.''

To make the exploration self-contained, we shall consider a slightly
generalized form of the equation system (\ref{1.3}) involving some
additional immigration structures. In Section \ref{sec2}, we present some
preliminary results on inhomogeneous immigration superprocesses and
CBI-processes. In Section \ref{sec3} a class of CBI-processes with predictable
immigration rates are constructed as pathwise unique solutions of
stochastic integral equations driven by time--space noises. In Section
\ref{sec4} we
introduce the path-valued increasing Markov processes and identify them as
path-valued branching processes with immigration. A construction of those
processes is given in Section \ref{sec5} using a system of stochastic equations
generalizing (\ref{1.3}). In Section \ref{sec6} we derive a homogeneous nonlocal
branching immigration superprocess from the flow. The properties of the
process under an excursion law are studied in Section~\ref{sec7}.

We sometimes write $\mbb{R}_+$ for $[0,\infty)$. Let $F(T)$ denote
the set
of positive right continuous increasing functions on an interval
$T\subset
\mbb{R}$. For a measure $\mu$ and a function $f$ on a measurable
space we
write $\langle\mu,f\rangle= \int fd\mu$ if the integral exists.
Throughout this
paper, we make the conventions
\[
\int_a^b = \int_{(a,b]}
\quad\mbox{and}\quad \int_a^\infty= \int_{(a,\infty)}
\]
for any $b\ge a\in\mbb{R}$. Other notations are explained as they first
appear.

%%%%%%%%%%%% (Section 2) %%%%%%%%%%%%%%

%s2 #&#
\section{Inhomogeneous immigration superprocesses}\label{sec2}

In this section, we present some preliminary results on inhomogeneous
immigration superprocesses and CBI-processes. Suppose that $T\subset
\mbb{R}$ is an interval, and $E$ is a Lusin topological space. Let
$\tilde{E} = T\times E$. A function $(s,x)\mapsto f(s,x)$ on $\tilde{E}$
is said to be \textit{locally bounded} if for each compact interval
$S\subset T$ the restriction of $(s,x)\mapsto f(s,x)$ to $S\times E$ is
bounded. Let $M(E)$ be the space of finite Borel measures on $E$ endowed
with the topology of weak convergence. Let $B^+(E)$ be the set of bounded
positive Borel functions on $E$. Let $\mcr{I}(E)$ denote the set of all
functionals $I$ on $B^+(E)$ with the representation
%
%e2.1 #&#
\begin{equation}
\label{2.1} I(f) = \langle\lambda,f\rangle+ \int_{M(E)^\circ}
\bigl(1-e^{-\langle\nu,f\rangle} \bigr) L(d\nu),\qquad f\in B^+(E),
\end{equation}
where $\lambda\in M(E)$ and $(1\land\langle\nu,1\rangle) L(d\nu)$
is a finite
measure on $M(E)^\circ:= M(E)\setminus\{0\}$. Let $\mcr{J}(E)$ denote
the set of all functionals on $B^+(E)$ of the form $f\mapsto J(f):=
a+I(f)$ with $a\ge0$ and $I\in\mcr{I}(E)$. By Theorems 1.35 and 1.37 in
\citet{Li11} one can prove the following:
%
%th2.1 #&#
\begin{theorem}\label{t2.1} There is a one-to-one correspondence between
functionals $V\in\mcr{J}(E)$ and infinitely divisible sub-probability
measures $Q$ on $M(E)$, which is determined by
%
%e2.2 #&#
\begin{equation}
\label{2.2} \int_{M(E)} e^{-\langle\nu,f\rangle} Q(d\nu) = \exp\bigl\{
-V(f) \bigr\},\qquad f\in B^+(E).
\end{equation}
\end{theorem}
%
%th2.2 #&#
\begin{theorem}\label{t2.2} If $U\in\mcr{J}(E)$ and if $V\dvtx  f\mapsto
v(\cdot,f)$
is an operator on $B^+(E)$ such that $v(x,\cdot)\in\mcr{J}(E)$ for all
$x\in E$, then $U\circ V\in\mcr{J}(E)$.
\end{theorem}

Suppose that $(P_{r,t}\dvtx  t\ge r\in T)$ is an inhomogeneous Borel right
transition semigroup on $E$. Let $\xi= (\Omega, \mcr{F}, \mcr{F}_{r,t},
\xi_t, \mbf{P}_{r,x})$ be a right continuous inhomogeneous Markov process
realizing $(P_{r,t}\dvtx  t\ge r\in T)$. Let $(s,x)\mapsto b_s(x)$ be a Borel
function on $\tilde{E}$, and let $(s,x)\mapsto c_s(x)$ be a positive Borel
function on $\tilde{E}$. Let $\eta_s(x,dy)$ be a kernel from $\tilde{E}$
to $E$, and let $H_s(x,d\nu)$ be a kernel from $\tilde{E}$ to
$M(E)^\circ$.
Suppose that the function
\[
\bigl|b_s(x)\bigr| + c_s(x) + \eta_s(x,E) + \int
_{M(E)^\circ} \bigl(\langle\nu,1\rangle\land\langle\nu,1
\rangle^2 + \langle\nu_x,1\rangle\bigr)
H_s(x,d\nu)
\]
on $S\times E$ is locally bounded, where $\nu_x(dy)$ denotes the
restriction of $\nu(dy)$ to $E\setminus\{x\}$. For $(s,x)\in\tilde{E}$
and $f\in B^+(E)$ define
%
%e2.3 #&#
\begin{eqnarray}
\label{2.3} \phi_s(x,f) &=& b_s(x)f(x) +
c_s(x)f(x)^2 - \int_E f(y)
\eta_s(x,dy)
\nonumber\\[-8pt]\\[-8pt]
&&{} + \int_{M(E)^\circ} \bigl[e^{-\langle\nu,f\rangle}-1+\nu\bigl(\{x\}
\bigr)f(x) \bigr] H_s(x,d\nu).
\nonumber
\end{eqnarray}
Let $T_t = T\cap(-\infty,t]$ for $t\in T$. By Theorem 6.10 in \citet{Li11}
one can show there is an inhomogeneous Borel right transition semigroup
$(Q_{r,t}\dvtx  t\ge r\in T)$ on the state space $M(E)$ defined by
%
%e2.4 #&#
\begin{equation}
\label{2.4} \int_{M(E)} e^{-\langle\nu,f\rangle} Q_{r,t}(\mu,d
\nu) = \exp\bigl\{-\langle\mu,V_{r,t}f\rangle\bigr\},\qquad f\in B^+(E),
\end{equation}
where $(r,x)\mapsto v_{r,t}(x):= V_{r,t}f(x)$ is the unique locally
bounded positive solution to the integral equation
%
%e2.5 #&#
\begin{equation}
\label{2.5} v_{r,t}(x) = \mbf{P}_{r,x}\bigl[f(
\xi_t)\bigr] - \int_r^t
\mbf{P}_{r,x} \bigl[\phi_s(\xi_s,v_{s,t})
\bigr]\,ds,\qquad r\in T_t, x\in E.
\end{equation}

Let us consider a right continuous realization $X = (W, \mcr{G},
\mcr{G}_{r,t}, X_t, \mbf{Q}_{r,\mu})$ of the transition semigroup
$(Q_{r,t}\dvtx  t\ge r\in T)$ defined by (\ref{2.4}). Suppose that
$(s,x)\mapsto g_s(x)$ is a locally bounded positive Borel function on
$\tilde{E}$. Let $\psi_s(x,f) = - g_s(x) + \phi_s(x,f)$ for $f\in B^+(E)$.
Following the proofs of Theorems 5.15 and 5.16 in \citet{Li11}, one can see
%
%e2.6 #&#
\begin{equation}
\label{2.6} \mbf{Q}_{r,\mu} \exp\biggl\{-\langle X_t,f
\rangle- \int_r^t \langle X_s,g_s
\rangle \,ds \biggr\} = \exp\bigl\{-\langle\mu,U_{r,t}f\rangle\bigr\},
\end{equation}
where $(r,x)\mapsto u_{r,t}(x):= U_{r,t}f(x)$ is the unique locally
bounded positive solution to
%
%e2.7 #&#
\begin{equation}
\label{2.7}\quad u_{r,t}(x) = \mbf{P}_{r,x}\bigl[f(
\xi_t)\bigr] - \int_r^t
\mbf{P}_{r,x} \bigl[\psi_s(\xi_s,u_{s,t})
\bigr]\,ds,\qquad r\in T_t, x\in E.
\end{equation}
Then there is an inhomogeneous Borel right sub-Markov transition semigroup
$(Q_{r,t}^g\dvtx  t\ge r\in T)$ on $M(E)$ given by
%
%e2.8 #&#
\begin{equation}
\label{2.8} \int_{M(E)} e^{-\langle\nu,f\rangle} Q_{r,t}^g(
\mu,d\nu) = \exp\bigl\{-\langle\mu,U_{r,t}f\rangle\bigr\}.
\end{equation}
A Markov process with transition semigroup given by (\ref{2.8}) is called
an \textit{inhomogeneous superprocess} with \textit{branching mechanisms}
$\{\psi_s\dvtx  s\in T\}$. The family of operators $(U_{r,t}\dvtx  t\ge r\in T)$ is
called the \textit{cumulant semigroup} of the superprocess. From
(\ref{2.8}) one can derive the following \textit{branching property}:
%
%e2.9 #&#
\begin{equation}
\label{2.9} Q_{r,t}^g(\mu_1+\mu_2,
\cdot) = Q_{r,t}^g(\mu_1,\cdot)*Q_{r,t}^g(
\mu_2,\cdot)
\end{equation}
for $t\ge r\in T$ and $\mu_1,\mu_2\in M(E)$, where ``$*$'' denotes the
convolution operation. Some special branching mechanisms are given in
\citet{DawGorLi02}, \citet{Dyn93} and
Li (\citeyear{Li92}, \citeyear{Li11}). Clearly, the
semigroup $(Q_{r,t}\dvtx  t\ge r\in T)$ given by (\ref{2.4}) corresponds to a
conservative inhomogeneous superprocess. In general, the inhomogeneous
superprocess is not necessarily conservative.

We can append an additional immigration structure to the inhomogeneous
superprocess. Suppose that $\rho(ds)$ is a Radon measure on $T$ and
$\{J_s\dvtx  s\in T\}\subset\mcr{J}(E)$ is a family of functionals such that
$s\mapsto J_s(f)$ is a locally bounded Borel function on $T$ for each
$f\in B^+(E)$.
%
%th2.3 #&#
\begin{theorem}\label{t2.3} There is an inhomogeneous transition semigroup
$(Q_{r,t}^{\rho,J}\dvtx\break  t\ge r\in T)$ on $M(E)$ given by
%
%e2.10 #&#
\begin{equation}
\label{2.10}\qquad\quad \int_{M(E)} e^{-\langle\nu,f\rangle} Q_{r,t}^{\rho,J}(
\mu,d\nu) = \exp\biggl\{-\langle\mu,U_{r,t}f\rangle- \int
_r^t J_s(U_{s,t}f)
\rho(ds) \biggr\},
\end{equation}
where $(r,x)\mapsto u_{r,t}(x):= U_{r,t}f(x)$ is the unique locally
bounded positive solution to (\ref{2.7}).
\end{theorem}
\begin{pf}
By Theorems \ref{t2.1} and \ref{t2.2}, for any $t\ge r\in T$
we can
define an infinitely divisible sub-probability measure $N_{r,t}$ on $M(E)$
by
\[
\int_{M(E)} e^{-\langle\nu,f\rangle} N_{r,t}(d\nu) = \exp
\biggl\{-\int_r^t J_s(U_{s,t}f)
\rho(ds) \biggr\}.
\]
It is easy to check that
\[
N_{r,t} = \bigl(N_{r,s}Q_{s,t}^g
\bigr)*N_{s,t},\qquad t\ge s\ge r\in T,
\]
where
\[
N_{r,s}Q_{s,t}^g = \int_{M(E)}
N_{r,s}(d\mu)Q_{s,t}^g(\mu,\cdot).
\]
Following the arguments in Li (\citeyear{Li01}, \citeyear{Li11}) one can show
%
%e2.11 #&#
\begin{equation}
\label{2.11} Q_{r,t}^{\rho,J}(\mu,\cdot) = Q_{r,t}^g(
\mu,\cdot)*N_{r,t},\qquad t\ge r\in T,
\end{equation}
defines an inhomogeneous sub-Markov transition semigroup on $M(E)$.
Clearly, the Laplace functional of this transition semigroup is given by
(\ref{2.10}).
\end{pf}

If a Markov process with state space $M(E)$ has transition semigroup
$(Q_{r,t}^{\rho,J}\dvtx  t\ge r\in T)$ given by (\ref{2.10}), we call it an
\textit{inhomogeneous immigration superprocess} with \textit{immigration
mechanisms} $\{J_s\dvtx  s\in T\}$ and \textit{immigration measure} $\rho
$. The
intuitive meaning of the model is clear in view of (\ref{2.11}). That is,
the population at any time $t\ge0$ is made up of two parts, the native
part generated by the mass $\mu\in M(E)$ at time $r\ge0$ has distribution
$Q_{r,t}^g(\mu,\cdot)$ and the immigration in the time interval $(r,t]$
gives the distribution $N_{r,t}$. When $E$ shrinks to a singleton, we can
identify $M(E)$ with the positive half line $\mbb{R}_+ = [0,\infty)$. In
this case, the transition semigroups given by (\ref{2.8}) and (\ref{2.10})
determine one-dimensional CB- and CBI-processes, respectively.

Now let us consider a branching mechanism $\phi$ of the form (\ref{1.1}).
We can define the transition semigroup $(P_t)_{t\ge0}$ of a homogeneous
CB-process by
%
%e2.12 #&#
\begin{equation}
\label{2.12} \int_{\mbb{R}_+} e^{-\lambda y} P_t(x,dy)
= e^{-xv_t(\lambda)},\qquad t,\lambda\ge0,
\end{equation}
where $t\mapsto v_t(\lambda)$ is the unique locally bounded positive
solution of
\[
v_t(\lambda) = \lambda- \int_0^t
\phi\bigl(v_s(\lambda)\bigr)\,ds,
\]
which is essentially a special form of (\ref{2.5}). We can write the above
integral equation into its differential form
%
%e2.13 #&#
\begin{equation}
\label{2.13} \frac{d}{dt}v_t(\lambda) = - \phi
\bigl(v_t(\lambda)\bigr),\qquad v_0(\lambda) = \lambda.
\end{equation}
The Chapman--Kolmogorov equation of $(P_t)_{t\ge0}$ implies
$v_r(v_t(\lambda)) = v_{r+t}(\lambda)$ for all $r,t,\lambda\ge0$.
The set
of functions $(v_t)_{t\ge0}$ is the \textit{cumulant semigroup}. Observe
that $\lambda\mapsto\phi(\lambda)$ is continuously differentiable with
\[
\phi'(\lambda) = b + \sigma^2\lambda+ \int
_0^\infty z\bigl(1-e^{-z\lambda}\bigr)m(dz),\qquad
\lambda\ge0.
\]
By differentiating (\ref{2.12}) and (\ref{2.13}) in $\lambda\ge0$
one can
show
%
%e2.14 #&#
\begin{equation}
\label{2.14} \int_{\mbb{R}_+} y P_t(x,dy) = x\,
\frac{d}{d\lambda} v_t(\lambda) \bigg|_{\lambda=0+} = xe^{-bt}.
\end{equation}

It is easy to see that $(P_t)_{t\ge0}$ is a Feller semigroup. Let us
consider a c\`adl\`ag realization $X = (\iittOmega, \mcr{F}, \mcr{F}_{r,t},
X_t, \mbf{P}_{r,x})$ of the corresponding CB-process with an arbitrary
initial time $r\ge0$. Let $\eta(ds)$ be a Radon measure on $[0,\infty)$.
By Theorem~5.15 in \citet{Li11}, for $t\ge r\ge0$ and $f\in B^+[0,t]$, we
have
%
%e2.15 #&#
\begin{equation}
\label{2.15} \mbf{P}_{r,x} \biggl[\exp\biggl\{-\int
_{[r,t]} f(s)X_s \eta(ds) \biggr\} \biggr] = \exp
\bigl\{-xu^t(r,f) \bigr\},
\end{equation}
where $r\mapsto u^t(r,f)$ is the unique bounded positive solution to
%
%e2.16 #&#
\begin{equation}
\label{2.16} u^t(r,f) + \int_r^t
\phi\bigl(u^t(s,f)\bigr)\,ds = \int_{[r,t]} f(s)
\eta(ds),\qquad 0\le r\le t.
\end{equation}
In particular, for $r\ge0$ and $f\in B^+[0,\infty)$ with compact support,
we have
%
%e2.17 #&#
\begin{equation}
\label{2.17} \mbf{P}_{r,x} \biggl[\exp\biggl\{-\int
_r^\infty f(s)X_s \,ds \biggr\} \biggr]
= \exp\bigl\{-xu(r,f) \bigr\},
\end{equation}
where $r\mapsto u(r,f)$ is the unique compactly supported bounded positive
function on $[0,\infty)$ solving
%
%e2.18 #&#
\begin{equation}
\label{2.18} u(r,f) + \int_r^\infty\phi
\bigl(u(s,f)\bigr)\,ds = \int_r^\infty f(s) \,ds,\qquad r
\ge0.
\end{equation}
It is not hard to see that $u(r,f)=0$ for $r>l_f:= \sup\{t\ge0\dvtx
f(t)>0\}$. For any $r\ge0$ let
\[
\sigma_r(X) = \int_r^\infty
X_s \,ds.
\]

%th2.4 #&#
\begin{theorem}\label{t2.5} Suppose that $\phi(\lambda)\to\infty$ as
$\lambda\to\infty$. Then for any $\lambda\ge0$ we have
%
%e2.19 #&#
\begin{equation}
\label{2.19} \mbf{P}_{r,x} \bigl[e^{-\lambda\sigma_r(X)}1_{\{\sigma
_r(X)< \infty
\}}
\bigr] = \exp\bigl\{-x\phi^{-1}(\lambda) \bigr\},
\end{equation}
where $\phi^{-1}$ is the right inverse of $\phi$ defined by
%
%e2.20 #&#
\begin{equation}
\label{2.20} \phi^{-1}(\lambda) = \inf\bigl\{z\ge0\dvtx  \phi(z)> \lambda
\bigr\}.
\end{equation}
\end{theorem}
\begin{pf}
A proof of (\ref{2.19}) was already given in \citet{AbrDel12}. We
here give a simple derivation of the result since the argument is also
useful to prove the next theorem. By (\ref{2.15}) and (\ref{2.16}), for
any $t\ge r$ and $z,\theta\ge0$ we have
\[
\mbf{P}_{r,x} \biggl[\exp\biggl\{-zX_t - \theta\int
_r^t X_s \,ds \biggr\} \biggr] =
\exp\bigl\{-xu^t(r,z,\theta) \bigr\},
\]
where $r\mapsto u^t(r,z,\theta)$ is the unique bounded positive solution
to
\[
u^t(r,z,\theta) + \int_r^t \phi
\bigl(u^t(s,z,\theta)\bigr)\,ds = z + \theta(t-r),\qquad 0\le r\le t.
\]
Then one can see $u^t(r,z,\phi(z)) = z$. It follows that
\[
\mbf{P}_{r,x} \biggl[\exp\biggl\{-zX_t - \phi(z)\int
_r^t X_s \,ds \biggr\} \biggr] =
e^{-zx}.
\]
Since $\sigma_r(X)< \infty$ implies $\lim_{t\to\infty}X_t = 0$, if
$\phi(z)>0$, we get
\[
\mbf{P}_{r,x} \bigl[e^{-\phi(z)\sigma_r(X)}1_{\{\sigma_r(X)< \infty
\}} \bigr] =
e^{-zx}.
\]
That gives (\ref{2.19}) first for $\lambda= \phi(z)> 0$ and then for all
$\lambda\ge0$.
\end{pf}

Let $t\mapsto\rho(t)$ be a locally bounded positive Borel function on
$[0,\infty)$. Suppose that $h\ge0$ is a constant and $zn(dz)$ is a finite
measure on $(0,\infty)$. Let $\psi$ be an \textit{immigration mechanism}
given by
%
%e2.21 #&#
\begin{equation}
\label{2.21} \psi(\lambda) = h\lambda+ \int_0^\infty
\bigl(1-e^{-z\lambda}\bigr)n(dz),\qquad \lambda\ge0.
\end{equation}
By Theorem \ref{t2.3} we can define an inhomogeneous transition semigroup
$\{P_{r,t}^\rho\dvtx\break   t\ge r\ge0\}$ on $\mbb{R}_+$ by
%
%e2.22 #&#
\begin{equation}
\label{2.22}\quad \int_{\mbb{R}_+} e^{-\lambda y} P_{r,t}^\rho(x,dy)
= \exp\biggl\{-xv_{t-r}(\lambda)-\int_r^t
\psi\bigl(v_{t-s}(\lambda)\bigr)\rho(s) \,ds \biggr\}.
\end{equation}
A positive Markov process with transition semigroup $(P_{r,t}^\rho
)_{t\ge
r\ge0}$ is called an inhomogeneous CBI-process with \textit{immigration
rate} $\rho= \{\rho(t)\dvtx  t\ge0\}$. It is easy to see that the homogeneous
time--space semigroup associated with $(P_{r,t}^\rho)_{t\ge r\ge0}$ is a
Feller transition semigroup. Then $(P_{r,t}^\rho)_{t\ge r\ge0}$ has a
c\`adl\`ag realization $Y = (\iittOmega, \mcr{F}, \mcr{F}_{r,t}, Y_t,
\mbf{P}_{r,x}^\rho)$. A modification of the proof of Theorem 5.15 in \citet{Li11} shows that, for $t\ge r\ge0$ and $f\in B^+[0,t]$,
%
%e2.23 #&#
\begin{eqnarray}
\label{2.23}
&&\mbf{P}_{r,x}^\rho\biggl[\exp\biggl\{- \int
_{[r,t]} f(s)Y_s \eta(ds) \biggr\} \biggr]
\nonumber\\[-8pt]\\[-8pt]
&&\qquad = \exp\biggl\{-xu^t(r,f) - \int_r^t
\psi\bigl(u^t(s,f)\bigr)\rho(s) \,ds \biggr\},
\nonumber
\end{eqnarray}
where $r\mapsto u^t(r,f)$ is the unique bounded positive solution to
(\ref{2.16}). In particular, for $r\ge0$ and $f\in B^+[0,\infty)$ with
compact support, we have
%
%e2.24 #&#
\begin{eqnarray}
\label{2.24}
&&\mbf{P}_{r,x}^\rho\biggl[\exp\biggl\{- \int
_r^\infty f(s)Y_s \,ds \biggr\} \biggr]
\nonumber\\[-8pt]\\[-8pt]
&&\qquad = \exp\biggl\{-xu(r,f) - \int_r^\infty\psi
\bigl(u(s,f)\bigr)\rho(s) \,ds \biggr\},
\nonumber
\end{eqnarray}
where $r\mapsto u(r,f)$ is the unique compactly supported bounded positive
solution to (\ref{2.18}). For any $r\ge0$ let
\[
\sigma_r(Y) = \int_r^\infty
Y_s \,ds.
\]
By a modification of the proof of Theorem \ref{t2.5}, we get the following:
%
%th2.5 #&#
\begin{theorem}\label{t2.6} Suppose that $\phi(\lambda)\to\infty$ as
$\lambda\to\infty$. Then for any $r,\lambda\ge0$ we have
\[
\mbf{P}_{r,x}^\rho\bigl[e^{-\lambda\sigma_r(Y)}1_{\{\sigma_r(Y)<
\infty\}}
\bigr] = \exp\biggl\{-x\phi^{-1}(\lambda) - \psi\bigl(
\phi^{-1}(\lambda)\bigr)\int_r^\infty
\rho(s) \,ds \biggr\},
\]
where $\phi^{-1}(\lambda)$ is defined by (\ref{2.20}).
\end{theorem}

%%%%%%%%%%%% (Section 3) %%%%%%%%%%%%%%

%s3 #&#
\section{The predictable immigration rate}\label{sec3}

The main purpose of this section is to give a construction of the
CBI-process with transition semigroup $(P_{r,t}^\rho)_{t\ge r\ge0}$
defined by (\ref{2.22}) as the pathwise unique solution of a
stochastic integral equation driven by time--space noises. For the
convenience of applications, we shall generalize the model slightly by
considering a random immigration rate. This is essential for our study
of the
path-valued Markov processes. The reader is referred to
\citet{BerLeG06}, Dawson and Li (\citeyear{DawLi06}, \citeyear{DawLi12}),
\citet{FuLi10} and \citet{LiMyt11} for some related results.

Suppose that $(\iittOmega, \mcr{F}, \mcr{F}_t, \mbf{P})$ is a filtered
probability space satisfying the usual hypotheses. Let $\{W(t,\cdot)\dvtx
t\ge
0\}$ be an $(\mcr{F}_t)$-white noise on $(0,\infty)$ based on the Lebesgue
measure and let $\{p_0(t)\dvtx  t\ge0\}$ and $\{p_1(t)\dvtx  t\ge0\}$ be
$(\mcr{F}_t)$-Poisson point processes on $(0,\infty)^2$ with
characteristic measures $m(dz)\,du$ and $n(dz)\,du$, respectively. We assume
that the white noise and the Poisson processes are independent of each
other. Let $W(ds,du)$ denote the stochastic integral on $(0,\infty)^2$
with respect to the white noise. Let $N_0(ds,dz,du)$ and $N_1(ds,dz,du)$
denote the Poisson random measures on $(0,\infty)^3$ associated with
$\{p_0(t)\}$ and $\{p_1(t)\}$, respectively. Let $\tilde{N}_0(ds,dz,du)$
denote the compensated random measure associated with $\{p_0(t)\}$.
Suppose that $\rho= \{\rho(t)\dvtx  t\ge0\}$ is a positive
$(\mcr{F}_t)$-predictable process such that $t\mapsto\mbf{P}[\rho
(t)]$ is
locally bounded. We are interested in positive c\`adl\`ag solutions of the
stochastic equation
%
%e3.1 #&#
\begin{eqnarray}
\label{3.1} Y_t &=& Y_0 + \sigma\int
_0^t\int_0^{Y_{s-}}
W(ds,du) + \int_0^t \int_0^\infty
\int_0^{Y_{s-}} z \tilde{N}_0(ds,dz,du)
\nonumber\\[-8pt]\\[-8pt]
&&{} + \int_0^t \bigl(h\rho(s)-bY_{s-}
\bigr) \,ds + \int_0^t \int_0^\infty
\int_0^{\rho(s)} z N_1(ds,dz,du).
\nonumber
\end{eqnarray}

For any positive c\`adl\`ag solution $\{Y_t\dvtx  t\ge0\}$ of (\ref{3.1})
satisfying $\mbf{P}[Y_0]< \infty$, one can use a standard stopping time
argument to show that $t\mapsto\mbf{P}[Y_t]$ is locally bounded and
%
%e3.2 #&#
\begin{equation}
\label{3.2} \mbf{P}[Y_t] = \mbf{P}[Y_0] +
\psi'(0)\int_0^t \mbf{P}\bigl[
\rho(s)\bigr] \,ds - b\int_0^t
\mbf{P}[Y_s] \,ds,
\end{equation}
where
\[
\psi'(0) = h + \int_0^\infty z
n(dz).
\]
By It\^o's formula, it is easy to see that $\{Y_t\dvtx  t\ge0\}$ solves the
following martingale problem: for every $f\in C^2(\mbb{R}_+)$,
%
%e3.3 #&#
\begin{eqnarray}
\label{3.3} f(Y_t) &=& f(Y_0) + \mbox{local mart.} - b
\int_0^t Y_sf'(Y_s)\,ds
+ \frac{1}{2} \sigma^2\int_0^t
Y_sf''(Y_s)\,ds
\nonumber
\\
&&{} + \int_0^tY_s\,ds\int
_0^\infty\bigl[f(Y_s + z) -
f(Y_s) - zf'(Y_s)\bigr]m(dz)
\\
&&{} + \int_0^t\rho(s) \biggl
\{hf'(Y_s) + \int_0^\infty
\bigl[f(Y_s + z) - f(Y_s)\bigr] n(dz) \biggr\}\,ds.
\nonumber
\end{eqnarray}

%pr3.1 #&#
\begin{prop}\label{t3.1} Suppose that $\{Y_t\dvtx  t\ge0\}$ is a positive
c\`adl\`ag solution of (\ref{3.1}) and $\{Z_t\dvtx  t\ge0\}$ is a positive
c\`adl\`ag solution of the equation with $(b,\rho)$ replaced by
$(c,\eta)$. Then we have
\begin{eqnarray*}
\mbf{P}\bigl[|Z_t-Y_t|\bigr] &\le& \mbf{P}\bigl[|Z_0-Y_0|\bigr]
+ \psi'(0) \int_0^t \mbf{P}
\bigl[\bigl|\eta(s)-\rho(s)\bigr|\bigr] \,ds
\\
&&{} + |c|\int_0^t \mbf{P}\bigl[|Z_s-Y_s|\bigr]
\,ds + |b-c|\int_0^t \mbf{P}[Y_s]
\,ds.
\end{eqnarray*}
\end{prop}
\begin{pf}
For each integer $n\ge0$ define $a_n = \exp\{-n(n+1)/2\}$.
Then $a_n\to0$ decreasingly as $n\to\infty$ and
\[
\int_{a_n}^{a_{n-1}}z^{-1} \,dz = n,\qquad n\ge1.
\]
Let $x\mapsto g_n(x)$ be a positive continuous function supported by
$(a_n,a_{n-1})$, so that
\[
\int_{a_n}^{a_{n-1}}g_n(x)\,dx=1
\]
and $g_n(x)\le2(nx)^{-1}$ for every $x>0$. Let
\[
f_n(z)=\int_0^{|z|}dy\int
_0^yg_n(x)\,dx,\qquad z\in\mbb{R}.
\]
It is easy to see that $|f_n'(z)|\le1$ and
\[
0\le|z|f_n''(z) = |z|g_n\bigl(|z|\bigr)
\le2n^{-1},\qquad z\in\mbb{R}.
\]
Moreover, we have $f_n(z)\rightarrow|z|$ increasingly as $n\to\infty$.
Let $\alpha_t = Z_t-Y_t$ for $t\ge0$. From (\ref{3.1}) we have
%
%e3.4 #&#
\begin{eqnarray}
\label{3.4} \alpha_t &=& \alpha_0 + h\int
_0^t \bigl[\eta(s)-\rho(s)\bigr] \,ds - c\int
_0^t \alpha_{s-} \,ds + (b-c)\int
_0^t Y_{s-} \,ds
\nonumber
\\
&&{} + \sigma\int_0^t\int
_{Y_{s-}}^{Z_{s-}} W(ds,du) + \int_0^t
\int_0^\infty\int_{Y_{s-}}^{Z_{s-}}
z \tilde{N}_0(ds,dz,du)
\\
&&{} + \int_0^t \int_0^\infty
\int_{\rho(s)}^{\eta(s)} z N_1(ds,dz,du).
\nonumber
\end{eqnarray}
By this and It\^o's formula,
%
%e3.5 #&#
\begin{eqnarray}
\label{3.5}\quad f_n(\alpha_t) &=& f_n(
\alpha_0) + h\int_0^t
f_n'(\alpha_s) \bigl[\eta(s)-\rho(s)
\bigr] \,ds - c\int_0^t f_n'(
\alpha_s)\alpha_s \,ds
\nonumber
\\
&&{} + (b-c)\int_0^t f_n'(
\alpha_s)Y_s \,ds + \frac{1}{2}\sigma^2
\int_0^t f_n''(
\alpha_s)|\alpha_s|\,ds
\nonumber
\\
&&{} + \int_0^t\alpha_s1_{\{\alpha_s>0\}}\,ds
\int_0^\infty\bigl[f_n(
\alpha_s+z) - f_n(\alpha_s) -
zf_n'(\alpha_s)\bigr] m(dz)
\nonumber
\\
&&{} - \int_0^t\alpha_s1_{\{\alpha_s<0\}}\,ds
\int_0^\infty\bigl[f_n(
\alpha_s-z) - f_n(\alpha_s) +
zf_n'(\alpha_s)\bigr] m(dz)
\\
&&{} + \int_0^t\bigl[\eta(s)-\rho(s)
\bigr]1_{\{\eta(s)>\rho(s)\}}\,ds \int_0^\infty
\bigl[f_n(\alpha_s+z) - f_n(
\alpha_s)\bigr] n(dz)
\nonumber
\\
&&{} - \int_0^t\bigl[\rho(s)-\eta(s)
\bigr]1_{\{\rho(s)>\eta(s)\}}\,ds \int_0^\infty
\bigl[f_n(\alpha_s-z) - f_n(
\alpha_s)\bigr] n(dz)
\nonumber
\\
&&{} + \mbox{mart.}\nonumber
\end{eqnarray}
It is easy to see that $|f_n(a+x) - f_n(a)|\le|x|$ for any
$a,x\in\mbb{R}$. If $ax\ge0$, we have
\[
|a|\bigl|f_n(a+x) - f_n(a) - xf_n'(a)\bigr|
\le\bigl(2|ax|\bigr)\land\bigl(n^{-1}|x|^2\bigr).
\]
Taking the expectation in both sides of (\ref{3.5}) gives
\begin{eqnarray*}
\mbf{P}\bigl[f_n(\alpha_t)\bigr] &\le& \mbf{P}
\bigl[f_n(\alpha_0)\bigr] + h\int_0^t
\mbf{P}\bigl[\bigl|\eta(s)-\rho(s)\bigr|\bigr] \,ds + |c|\int_0^t
\mbf{P}\bigl[|\alpha_s|\bigr] \,ds
\\
&&{} + |b-c|\int_0^t \mbf{P}[Y_s]
\,ds + \int_0^t \mbf{P}\bigl[\bigl|\eta(s)-\rho(s)\bigr|
\bigr] \,ds\int_0^\infty z n(dz)
\\
&&{} + n^{-1}\sigma^2t + \int_0^t\,ds
\int_0^\infty\bigl\{\bigl(2z\mbf{P}\bigl[|
\alpha_s|\bigr]\bigr)\land\bigl(n^{-1}z^2\bigr)\bigr
\} m(dz).
\end{eqnarray*}
Then we get the desired estimate by letting $n\to\infty$.
\end{pf}
%
%pr3.2 #&#
\begin{prop}\label{t3.2} Suppose that $\{Y_t\dvtx  t\ge0\}$ is a positive
c\`adl\`ag solution of (\ref{3.1}), and $\{Z_t\dvtx  t\ge0\}$ is a positive
c\`adl\`ag solution of the equation with $(b,\rho)$ replaced by
$(c,\eta)$. Then we have
\begin{eqnarray*}
\mbf{P} \Bigl[\sup_{0\le s\le t}|Z_s-Y_s| \Bigr] &
\le& \mbf{P}\bigl[|Z_0-Y_0|\bigr] + \psi'(0) \int
_0^t \mbf{P}\bigl[\bigl|\eta(s)-\rho(s)\bigr|\bigr] \,ds
\\
&&{} + \biggl(|c|+2\int_1^\infty z m(dz) \biggr)
\int_0^t \mbf{P}\bigl[|Z_s-Y_s|\bigr]
\,ds
\\
&&{} + |b-c|\int_0^t \mbf{P}[Y_s]
\,ds + 2\sigma\biggl(\int_0^t
\mbf{P}\bigl[|Z_s-Y_s|\bigr]\,ds \biggr)^{1/2}
\\
&&{} + 2 \biggl(\int_0^t
\mbf{P}\bigl[|Z_s-Y_s|\bigr]\,ds \int_0^1
z^2 m(dz) \biggr)^{1/2}.
\end{eqnarray*}
\end{prop}
\begin{pf}
This follows by applying Doob's martingale inequality to
(\ref{3.4}).
\end{pf}
%
%th3.3 #&#
\begin{theorem}\label{t3.3} For any $Y_0\ge0$ there is a pathwise unique
positive c\`adl\`ag solution $\{Y_t\dvtx  t\ge0\}$ of (\ref{3.1}).
\end{theorem}
\begin{pf}
The pathwise uniqueness of the solution follows by
Proposition \ref{t3.1} and Gronwall's inequality. Without loss of
generality, we may assume $Y_0\ge0$ is deterministic in proving the
existence of the solution. We give the proof in three steps.

\textit{Step} 1. Let $B(t) = W((0,t]\times(0,1])$. Then $\{B(t)\dvtx
t\ge0\}$ is a standard Brownian motion. By Theorems 5.1 and 5.2 in \citet{DawLi06}, for any constant $\rho\ge0$ there is a pathwise unique
positive solution to
\begin{eqnarray*}
Y_t &=& Y_0 + \sigma\int_0^t
\sqrt{Y_{s-}} \,dB(s) + \int_0^t \int
_0^\infty\int_0^{Y_{s-}}
z \tilde{N}_0(ds,dz,du)
\\
&&{} + \int_0^t (h\rho-bY_{s-}) \,ds
+ \int_0^t \int_0^\infty
\int_0^{\rho} z N_1(ds,dz,du).
\end{eqnarray*}
It is simple to see that $\{Y_t\dvtx  t\ge0\}$ is a weak solution to
%
%e3.6 #&#
\begin{eqnarray}
\label{3.6} Y_t &=& Y_0 + \sigma\int
_0^t\int_0^{Y_{s-}}
W(ds,du) + \int_0^t \int_0^\infty
\int_0^{Y_{s-}} z \tilde{N}_0(ds,dz,du)
\nonumber\\[-8pt]\\[-8pt]
&&{} + \int_0^t (h\rho-bY_{s-}) \,ds
+ \int_0^t \int_0^\infty
\int_0^{\rho} z N_1(ds,dz,du).
\nonumber
\end{eqnarray}
As pointed out at the beginning of this proof, the pathwise uniqueness
holds for~(\ref{3.6}).

\textit{Step} 2. Let $0=r_0<r_1<r_2<\cdots$ be an increasing
sequence. For each $i\ge1$ let $\eta_i$ be a positive integrable random
variable measurable with respect to $\mcr{F}_{r_{i-1}}$. Let $\rho=
\{\rho(t)\dvtx  t\ge0\}$ be the positive $(\mcr{F}_t)$-predictable step
process given by
\[
\rho(t) = \sum_{i=1}^\infty
\eta_i 1_{(r_{i-1},r_i]}(t),\qquad t\ge0.
\]
By the result in the first step, on each interval $(r_{i-1},r_i]$ there is
a pathwise unique solution $\{Y_t\dvtx  r_{i-1}<t\le r_i\}$ to
\begin{eqnarray*}
Y_t &=& Y_{r_{i-1}} + \sigma\int_{r_{i-1}}^t
\int_0^{Y_{s-}} W(ds,du) + \int_{r_{i-1}}^t
\int_0^\infty\int_0^{Y_{s-}}
z \tilde{N}_0(ds,dz,du)
\\
&&{} + \int_{r_{i-1}}^t (h\eta_i-bY_{s-})
\,ds + \int_{r_{i-1}}^t \int_0^\infty
\int_0^{\eta_i} z N_1(ds,dz,du).
\end{eqnarray*}
Then $\{Y_t\dvtx  t\ge0\}$ is a solution to (\ref{3.1}).

\textit{Step} 3. Suppose that $\rho= \{\rho(t)\dvtx  t\ge0\}$ is
general positive $(\mcr{F}_t)$-predictable process such that $t\mapsto
\mbf{P}[\rho(t)]$ is locally bounded. Take a sequence of positive
predictable step processes $\rho_k = \{\rho_k(t)\dvtx  t\ge0\}$ so that
%
%e3.7 #&#
\begin{equation}
\label{3.7} \mbf{P} \biggl[\int_0^t \bigl|
\rho_k(s) - \rho(s)\bigr| \,ds \biggr] \to0
\end{equation}
for every $t\ge0$ as $k\to\infty$. Let $\{Y_k(t)\dvtx  t\ge0\}$ be the
solution to (\ref{3.1}) with $\rho=\rho_k$. By Proposition \ref{t3.1},
Gronwall's inequality and (\ref{3.7}), one sees
\[
\sup_{0\le s\le t}\mbf{P}\bigl[\bigl|Y_k(s)-Y_i(s)\bigr|\bigr]
\to0
\]
for every $t\ge0$ as $i,k\to\infty$. Then Proposition \ref{t3.2} implies
\[
\mbf{P} \Bigl[\sup_{0\le s\le t}\bigl|Y_k(s)-Y_i(s)\bigr|
\Bigr]\to0
\]
for every $t\ge0$ as $i,k\to\infty$. Thus there is a subsequence
$\{k_i\}\subset\{k\}$ and a c\`adl\`ag process $\{Y_t\dvtx  t\ge0\}$ so that
\[
\sup_{0\le s\le t}\bigl|Y_{k_i}(s)-Y_s\bigr|\to0
\]
almost surely for every $t\ge0$ as $i\to\infty$. It is routine to show
that $\{Y_t\dvtx  t\ge0\}$ is a solution to (\ref{3.1}).
\end{pf}
%
%th3.4 #&#
\begin{theorem}\label{t3.4} If $\rho= \{\rho(t)\dvtx  t\ge0\}$ is a deterministic
locally bounded positive Borel function, the solution $\{Y_t\dvtx  t\ge0\}$ of
(\ref{3.1}) is an inhomogeneous CBI-process with transition semigroup
$\{P_{r,t}^\rho\dvtx  t\ge r\ge0\}$ defined by (\ref{2.22}).
\end{theorem}
\begin{pf}
By the martingale problem (\ref{3.3}), when $\rho(t) = \rho$
is a
deterministic constant function, the process $\{Y_t\dvtx  t\ge0\}$ is a Markov
process with transition semigroup $\{P_{r,t}^\rho\dvtx  t\ge r\ge0\}$; see,
for example, Theorem 9.30 in \citet{Li11}. If $\rho= \{\rho(t)\dvtx  t\ge0\}
$ is a
general deterministic locally bounded positive Borel function, we can
take each
step function $\rho_k = \{\rho_k(t)\dvtx  t\ge0\}$ in the last proof to be
deterministic. Then the solution $\{Y_k(t)\dvtx  t\ge0\}$ of (\ref{3.1}) with
$\rho=\rho_k$ is an inhomogeneous CBI-process with transition semigroup
$\{P_{r,t}^{\rho_k}\dvtx  t\ge r\ge0\}$. In other words, for any $\lambda
\ge
0$, $t\ge r\ge0$ and $G\in\mcr{F}_r$ we have
\[
\mbf{P}\bigl[1_Ge^{-\lambda Y_k(t)}\bigr] = \mbf{P}
\biggl[1_G\exp\biggl\{-Y_k(r)v_{t-r}(\lambda)
- \int_r^t\rho_k(s)\psi
\bigl(v_{t-s}(\lambda)\bigr)\,ds \biggr\} \biggr].
\]
Letting $k\to\infty$ along the sequence $\{k_i\}$ mentioned in the last
proof gives
\[
\mbf{P}\bigl[1_Ge^{-\lambda Y_t}\bigr] = \mbf{P}
\biggl[1_G\exp\biggl\{-Y_rv_{t-r}(\lambda) -
\int_r^t\rho(s)\psi\bigl(v_{t-s}(
\lambda)\bigr)\,ds \biggr\} \biggr].
\]
Then $\{Y_t\dvtx  t\ge0\}$ is a CBI-process with immigration rate $\rho=
\{\rho(t)\dvtx  t\ge0\}$.
\end{pf}

In view of the result of Theorem \ref{t3.4}, the solution $\{Y_t\dvtx  t\ge
0\}$ to (\ref{3.1}) can be called an inhomogeneous CBI-process with
branching mechanism $\phi$, immigration mechanism $\psi$ and
\textit{predictable immigration rate} \mbox{$\rho= \{\rho(t)\dvtx  t\ge0\}$}.

%%%%%%%%%%%% (Section 4) %%%%%%%%%%%%%%

%s4 #&#
\section{Path-valued branching processes}\label{sec4}

In this section, we introduce some path-valued Markov processes, which are
essentially special forms of the immigration superprocesses defined by
(\ref{2.7}) and (\ref{2.10}). Suppose that $T\subset\mbb{R}$ is an
interval, and $\{\phi_q\dvtx  q\in T\}$ is a family of branching mechanisms,
where $\phi_q$ is given by (\ref{1.1}) with the parameters $(b,m) =
(b_q,m_q)$ depending on $q\in T$. We call $\{\phi_q\dvtx  q\in T\}$ an
\textit{admissible family} if for each $\lambda\ge0$, the function
$q\mapsto\phi_q(\lambda)$ is decreasing and continuously differentiable
with the derivative $\zeta_q(\lambda):= - (\partial/\partial
q)\phi_q(\lambda)$ of the form
%
%e4.1 #&#
\begin{equation}
\label{4.1} \zeta_q(\lambda) = \beta_q\lambda+ \int
_0^\infty\bigl(1-e^{-z\lambda
}\bigr)
n_q(dz),\qquad q\in T, \lambda\ge0,
\end{equation}
where $\beta_q\ge0$ and $n_q(dz)$ is a $\sigma$-finite kernel from
$T$ to
$(0,\infty)$ satisfying
\[
\sup_{p\le\theta\le q} \biggl[\beta_\theta+ \int_0^\infty
z n_\theta(dz) \biggr]< \infty,\qquad q\ge p\in T.
\]
For an admissible family $\{\phi_q\dvtx  q\in T\}$, we clearly have
%
%e4.2 #&#
\begin{equation}
\label{4.2} \phi_{p,q}(\lambda):= \phi_p(\lambda) -
\phi_q(\lambda) = \int_p^q
\zeta_\theta(\lambda)\,d\theta.
\end{equation}
It follows that
%
%e4.3 #&#
\begin{equation}
\label{4.3} b_q = b_p - \int_p^q
\beta_\theta \,d\theta- \int_p^q d
\theta\int_0^\infty z n_\theta(dz)
\end{equation}
and
%
%e4.4 #&#
\begin{equation}
\label{4.4} m_q(dz) = m_p(dz) + \int
_{\{p<\theta\le q\}} n_\theta(dz)\,d\theta.
\end{equation}
We say $q_0\in T$ is a \textit{critical point} of the admissible family
$\{\phi_q\dvtx  q\in T\}$ if $b_{q_0} = 0$, which means $\phi_{q_0}$ is a
critical branching mechanism. By (\ref{4.3}) one can see $q\mapsto
b_q$ is
a continuous decreasing function on $T$, so the set of critical points
$T_0\subset T$ can only be an interval.

Let us consider a function $\mu\in F(T)$ and an admissible family of
branching mechanisms $\{\phi_q\dvtx  q\in T\}$. Write $\mu(p,q] = \mu(q) -
\mu(p)$ for $q\ge p\in T$. Recall that (\ref{2.22}) defines the
transition semigroup $\{P_{r,t}^\rho\dvtx  t\ge r\ge0\}$ of an inhomogeneous
CBI-process $\{Y_t\dvtx  t\ge0\}$. Let $\mbf{P}_x^\rho(\phi,\psi,dw)$ denote
the law on $D^+[0,\infty)$ of such a process with initial value
$Y_0=x\ge
0$. Given any $\rho\in D^+[0,\infty)$, we define the probability measure
$\mbf{P}_{p,q}(\rho,dw)$ on $D^+[0,\infty)$ by
%
%e4.5 #&#
\begin{equation}
\label{4.5} \mbf{P}_{p,q}(\rho,B) = \int_{D^+[0,\infty)}
1_B(\rho+w) \mbf{P}_{\mu(p,q]}^\rho(
\phi_q,\phi_{p,q},dw)
\end{equation}
for Borel sets $B\subset D^+[0,\infty)$. In view of (\ref{2.24}), for any
$f\in B^+[0,\infty)$ with compact support, we have
%
%e4.6 #&#
\begin{eqnarray}
\label{4.6} &&\int_{D^+[0,\infty)} \exp\biggl\{-\int
_0^\infty f(s)w(s)\,ds \biggr\} \mbf{P}_{p,q}(
\rho,dw)
\nonumber\\[-8pt]\\[-8pt]
&&\qquad = \exp\biggl\{-\mu(p,q]u_q(0,f) - \int_0^\infty
u_{p,q}(s,f)\rho(s) \,ds \biggr\},
\nonumber
\end{eqnarray}
where $s\mapsto u_q(s):= u_q(s,f)$ is the unique compactly supported
bounded positive solution to
%
%e4.7 #&#
\begin{equation}
\label{4.7} u_q(s) + \int_s^\infty
\phi_q\bigl(u_q(t)\bigr)\,dt = \int_s^\infty
f(t) \,dt
\end{equation}
and
%
%e4.8 #&#
\begin{equation}
\label{4.8} u_{p,q}(s,f) = f(s) + \phi_{p,q}
\bigl(u_q(s,f)\bigr),\qquad s\ge0.
\end{equation}
We remark that $u_q(s,f) = u_{p,q}(s,f) =0$ for $s>l_f:= \sup\{t\ge0\dvtx
f(t)>0\}$.
%
%pr4.1 #&#
\begin{prop}\label{t4.1} For any $f\in B^+[0,\infty)$ with compact
support, we have
%
%e4.9 #&#
\begin{equation}
\label{4.9} u_p\bigl(s,u_{p,q}(\cdot,f)\bigr) =
u_q(s,f),\qquad s\ge0, p\le q\in T,
\end{equation}
and
%
%e4.10 #&#
\begin{equation}
\label{4.10} u_{p,\theta}\bigl(s,u_{\theta,q}(\cdot,f)\bigr) =
u_{p,q}(s,f),\qquad s\ge0, p\le\theta\le q\in T.
\end{equation}
\end{prop}
\begin{pf}
From (\ref{4.2}) and (\ref{4.7}) we can see that $s\mapsto
v(s):=
u_q(s,f)$ is a solution of
%
%e4.11 #&#
\begin{equation}
\label{4.11} v(s) = \int_s^\infty\bigl[f(t) +
\phi_{p,q}\bigl(u_q(t,f)\bigr)\bigr] \,dt - \int
_s^\infty\phi_p\bigl(v(t)\bigr) \,dt.
\end{equation}
On the other hand, by (\ref{4.7}) and (\ref{4.8}) we have
\begin{eqnarray*}
u_p\bigl(s,u_{p,q}(\cdot,f)\bigr) &=& \int
_s^\infty u_{p,q}(t,f)\,dt - \int
_s^\infty\phi_p\bigl(u_p
\bigl(t,u_{p,q}(\cdot,f)\bigr)\bigr) \,dt
\\
&=& \int_s^\infty\bigl[f(t)+\phi_{p,q}
\bigl(u_q(t,f)\bigr)\bigr] \,dt
\\
&&{} - \int_s^\infty\phi_p
\bigl(u_p\bigl(t,u_{p,q}(\cdot,f)\bigr)\bigr) \,dt.
\end{eqnarray*}
Then $s\mapsto u_p(s,u_{p,q}(\cdot,f))$ is also a solution to
(\ref{4.11}). By the uniqueness of the solution to the equation, we get
(\ref{4.9}). It follows that
\begin{eqnarray*}
u_{p,\theta}\bigl(s,u_{\theta,q}(\cdot,f)\bigr) &=&
u_{\theta,q}(s,f) + \phi_{p,\theta}\bigl(u_\theta
\bigl(s,u_{\theta,q}(\cdot,f)\bigr)\bigr)
\\
&=& f(s) + \phi_{\theta,q}\bigl(u_q(s,f)\bigr) +
\phi_{p,\theta}\bigl(u_q(s,f)\bigr)
\\
&=& f(s) + \phi_{p,q}\bigl(u_q(s,f)\bigr).
\end{eqnarray*}
Then we have (\ref{4.10}).
\end{pf}
%
%pr4.2 #&#
\begin{prop}\label{t4.2} For any $f\in B^+[0,\infty)$ with compact
support we have
%
%e4.12 #&#
\begin{equation}
\label{4.12}\qquad u_{p,q}(s,f) = f(s) + \int_p^q
\psi_\theta\bigl(s,u_{\theta,q}(\cdot,f)\bigr) \,d\theta,\qquad s\ge0, q\ge
p\in T,
\end{equation}
where $\psi_\theta(s,f) = \zeta_\theta(u_\theta(s,f))$.
\end{prop}
\begin{pf}
By (\ref{4.2}) and (\ref{4.8}) one can see $p\mapsto u_{p,q}(s,f)$
is a decreasing function. In view of (\ref{4.9}) and (\ref{4.10}), for
$q>\theta> p\in T$, we get
\[
u_{p,q}(s,f) = u_{p,\theta}\bigl(s,u_{\theta,q}(\cdot,f)
\bigr) = u_{\theta,q}(s,f)+\phi_{p,\theta}\bigl(u_\theta
\bigl(s,u_{\theta,q}(\cdot,f)\bigr)\bigr).
\]
Then we differentiate both sides to see
\[
\frac{d}{dp}u_{p,q}(s,f) \bigg|_{p=\theta-} = \frac{d}{dp}
\phi_{p,\theta} \bigl(u_\theta\bigl(s,u_{\theta,q}(\cdot,f)
\bigr)\bigr) \bigg|_{p=\theta-} = -\zeta_\theta\bigl(u_\theta
\bigl(s,u_{\theta,q}(\cdot,f)\bigr)\bigr),
\]
which implies (\ref{4.12}).\vadjust{\goodbreak}
\end{pf}

From (\ref{4.6}) one can see that $\mbf{P}_{p,q}(\rho,dw)$ is a
probability kernel on $D^+[0,\infty)$. By (\ref{4.9}) and (\ref
{4.10}) it
is easy to check that the family of kernels $\{\mbf{P}_{p,q}\dvtx  q\ge
p\in
T\}$ satisfies the Chapman--Kolmogorov equation. Then $\{\mbf{P}_{p,q}\dvtx
q\ge p\in T\}$ form an inhomogeneous Markov transition semigroup on
$D^+[0,\infty)$. This semigroup is closely related to some nonlocal
branching superprocesses. For $\alpha\ge0$ let $M[0,\alpha]$ be the
space of finite Borel measures on $[0,\alpha]$ furnished with the
topology of weak convergence.
%
%th4.3 #&#
\begin{theorem}\label{t4.3} There is a Markov transition semigroup
$\{\mbf{Q}_{p,q}^\alpha\dvtx  q\ge p\in T\}$ on $M[0,\alpha]$ such that, for
$f\in B^+[0,\alpha]$,
%
%e4.13 #&#
\begin{equation}
\label{4.13}\qquad \int_{M[0,\alpha]} e^{-\langle\nu,f\rangle} \mbf
{Q}_{p,q}^\alpha(\eta,d\nu) = \exp\bigl\{-
\mu(p,q]u_q^\alpha(0,f) - \bigl\langle\eta,u_{p,q}^\alpha(
\cdot,f)\bigr\rangle\bigr\},
\end{equation}
where
%
%e4.14 #&#
\begin{equation}
\label{4.14} u_q^\alpha(s,f) = u_q(s,f1_{[0,\alpha]}),\qquad
u_{p,q}^\alpha(s,f) = u_{p,q}(s,f1_{[0,\alpha]}).
\end{equation}
\end{theorem}
\begin{pf}
We first consider an absolutely continuous measure $\eta\in
M[0,\alpha]$ with a density $\rho\in D^+[0,\alpha]$. Suppose that $\{X_t\dvtx
t\ge0\}$ is a random path with distribution $\mbf{P}_{p,q}(\rho
1_{[0,\alpha]},\cdot)$ on $D^+[0,\infty)$. Let $\mbf
{Q}_{p,q}^\alpha
(\eta,\cdot)$ be the distribution on $M[0,\alpha]$ of the random measure
$X$ such that $X(dt) = X_t\,dt$ for $0\le t\le\alpha$. The Laplace
function of $\mbf{Q}_{p,q}^\alpha(\eta,\cdot)$ is clearly given by
(\ref{4.13}) and (\ref{4.14}). In particular, we can use those two
formulas to define a probability measure on $M[0,\alpha]$. For an
arbitrary $\eta\in M[0,\alpha]$, choose a sequence of absolutely
continuous measures $\{\eta_n\}\subset M[0,\alpha]$ with densities in
$D^+[0,\alpha]$ so that $\eta_n\to\eta$ weakly. Let
$\mbf{Q}_{p,q}^\alpha(\eta_n,\cdot)$ be the probability measure on
$M[0,\alpha]$ defined by
\[
\int_{M[0,\alpha]} e^{-\langle\nu,f\rangle} \mbf{Q}_{p,q}^\alpha
(\eta_n,d\nu) = \exp\bigl\{-\mu(p,q]u_q^\alpha(0,f)
- \bigl\langle\eta_n,u_{p,q}^\alpha(\cdot,f)\bigr
\rangle\bigr\}.
\]
For $f\in C^+[0,\alpha]$ one can see from (\ref{4.7}) and (\ref{4.8})
that $u_{p,q}^\alpha(\cdot,f)\in C^+[0,\alpha]$, and hence
\[
\lim_{n\to\infty}\int_{M[0,\alpha]} e^{-\langle\nu,f\rangle}
\mbf{Q}_{p,q}^\alpha(\eta_n,d\nu) = \exp\bigl\{-
\mu(p,q]u_q^\alpha(0,f) - \bigl\langle\eta,u_{p,q}^\alpha(
\cdot,f)\bigr\rangle\bigr\}.
\]
Then (\ref{4.13}) really gives the Laplace functional of a probability
measure $\mbf{Q}_{p,q}^\alpha(\eta,\cdot)$ on $M[0,\alpha]$ which
is the
weak limit of $\mbf{Q}_{p,q}^\alpha(\eta_n,\cdot)$ as $n\to\infty
$. It
is\vspace*{1pt} easy to see that $\mbf{Q}_{p,q}^\alpha(\eta,d\nu)$ is a kernel on
$M[0,\alpha]$. The semigroup property of the family
$\{\mbf{Q}_{p,q}^\alpha\dvtx  q\ge p\in T\}$ follows from (\ref{4.9}) and
(\ref{4.10}).
\end{pf}
%
%th4.4 #&#
\begin{theorem}\label{t4.4} Let $q\in T$ and $f\in B^+[0,\alpha]$. Then
$(p,s)\mapsto u_{p,q}^\alpha(s):= u_{p,q}^\alpha(s,f)$ is the unique,
locally bounded positive solution to
%
%e4.15 #&#
\begin{equation}
\label{4.15} u_{p,q}^\alpha(s) = f(s) + \int
_p^q \psi_\theta^\alpha
\bigl(s, u_{\theta,q}^\alpha\bigr) \,d\theta,\qquad s\in[0,\alpha], q\ge p\in
T,
\end{equation}
where $\psi_\theta^\alpha(s,f) = \zeta_\theta(u_\theta^\alpha(s,f))$.
Moreover, the transition semigroup $\{\mbf{Q}_{p,q}^\alpha\dvtx  q\ge p\in
T\}$ defines an immigration superprocess in $M[0,\alpha]$ with branching
mechanisms $\{-\psi_\theta^\alpha\dvtx  \theta\in T\}$, immigration mechanisms
$\{u_\theta^\alpha(0,\cdot)\dvtx  \theta\in T\}$ and immigration measure
$\mu$.
\end{theorem}
\begin{pf}
From (\ref{4.12}) one can see that $u_{p,q}^\alpha(s) =
u_{p,q}^\alpha(s,f)$ satisfies (\ref{4.15}). By letting $t = \alpha$ and
$\eta(ds) = ds$ in (\ref{2.15}), we infer that the functional
$f\mapsto
u_\theta^\alpha(s,f)$ on $B^+[0,\alpha]$ is the Laplace exponent of an
infinitely divisible probability measure carried by $M[s,\alpha]$. It is
easy to see that $\psi_\theta^\alpha(s,0) = 0$. By Theorem \ref
{t2.2} the
composed functional $f\mapsto\psi_\theta^\alpha(s,f)$ is also the
Laplace exponent of an infinitely divisible probability measure on
$M[s,\alpha]$. Then it has the representation
%
%e4.16 #&#
\begin{equation}
\label{4.16} \psi_\theta^\alpha(s,f) = \bigl\langle
\eta_\theta^\alpha(s),f\bigr\rangle+ \int_{M[s,\alpha]^\circ}
\bigl(1-e^{-\langle\nu,f\rangle}\bigr) H_\theta^\alpha(s,d\nu),
\end{equation}
where $\eta_\theta^\alpha(s)\in M[s,\alpha]$ and $(1\land\langle
\nu,1\rangle)
H_\theta^\alpha(s,d\nu)$ is a finite measure on $M[s,\alpha]^\circ
$. By
letting $f(t)=\lambda$ and taking the derivatives in both sides of
(\ref{4.16}), we have
\[
\frac{d}{d\lambda} \psi_\theta^\alpha(s,\lambda)
\bigg|_{\lambda=0+} = \bigl\langle\eta_\theta^\alpha(s),1\bigr
\rangle+ \int_{M[s,\alpha]^\circ
} \langle\nu,1\rangle H_\theta^\alpha(s,d
\nu).
\]
On the other hand, using (\ref{2.14}) and (\ref{2.15}),
\[
\frac{d}{d\lambda} u_\theta^\alpha(s,\lambda) \bigg|_{\lambda=0+}
= \int_s^\alpha e^{-b_\theta(t-s)} \,dt.
\]
From (\ref{4.1}) we have
\[
\frac{d}{d\lambda} \zeta_\theta(\lambda) \bigg|_{\lambda=0+} =
\beta_\theta+ \int_0^\infty z
n_\theta(dz).
\]
It follows that
\[
\frac{d}{d\lambda} \psi_\theta^\alpha(s,\lambda)
\bigg|_{\lambda=0+} = \biggl[\beta_\theta+ \int_0^\infty
z n_\theta(dz) \biggr]\int_s^\alpha
e^{-b_\theta(t-s)} \,dt.
\]
As a function of $(\theta,s)$, the above quantity is bounded on
$S\times
[0,\alpha]$ for each bounded closed interval $S\subset T$. By Example 2.5
of \citet{Li11} one sees that $f\mapsto-\psi_\theta^\alpha(\cdot,f)$
is a
special form of the operator given by (\ref{2.3}), and so (\ref
{4.15}) is
a special case of (\ref{2.7}). Thus $(p,s)\mapsto u_{p,q}^\alpha
(s,f)$ is
the unique locally bounded positive solution to (\ref{4.15}). By
(\ref{4.9}) we have
%
%e4.17 #&#
\begin{equation}
\label{4.17}\qquad \mu(p,q]u_q^\alpha(0,f) = \int
_p^q u_q^\alpha(0,f) \mu(d
\theta) = \int_p^q u_\theta^\alpha
\bigl(0,u_{\theta,q}^\alpha(\cdot,f)\bigr) \mu(d\theta).
\end{equation}
Then $\{\mbf{Q}_{p,q}^\alpha\dvtx  q\ge p\in T\}$ defines an immigration
superprocess in $M[0,\alpha]$ with branching mechanisms
$\{-\psi_\theta^\alpha\dvtx  \theta\in T\}$, immigration mechanisms
$\{u_\theta^\alpha(0,\cdot)\dvtx  \theta\in T\}$ and immigration measure
$\mu$.
\end{pf}

Let $\mcr{M}[0,\infty)$ denote the space of Radon measures on
$[0,\infty)$ endowed with the topology of vague convergence. For any
$\alpha\ge0$ we regard $M[0,\alpha]$ as the subset of
$\mcr{M}[0,\infty)$ consisting of the measures supported by
$[0,\alpha]$.
We can also embed $D^+[0,\infty)$ continuously into $\mcr{M}[0,\infty)$
by identifying the path $w\in D^+[0,\infty)$ and the measure $\nu\in
\mcr{M}[0,\infty)$ such that $\nu(ds) = w(s)\,ds$ for $s\ge0$.
%
%th4.5 #&#
\begin{theorem}\label{t4.5} There is an extension $\{\mbf{Q}_{p,q}\dvtx
q\ge
p\in
T\}$ of $\{\mbf{P}_{p,q}\dvtx  q\ge p\in T\}$ on $\mcr{M}[0,\infty)$,
which is
given by
%
%e4.18 #&#
\begin{equation}
\label{4.18}\hspace*{28pt} \int_{\mcr{M}[0,\infty)} e^{-\langle\nu,f\rangle} \mbf
{Q}_{p,q}(\eta,d\nu) = \exp\bigl\{-\mu(p,q]u_q(0,f) -
\bigl\langle\eta,u_{p,q}(\cdot,f)\bigr\rangle\bigr\}
\end{equation}
for $f\in B^+[0,\infty)$ with compact support.
\end{theorem}
\begin{pf}
Given $\eta\in\mcr{M}[0,\infty)$, we define $\pi_\alpha\eta\in
M[0,\alpha]$ by $\pi_\alpha\eta(ds) =\break  1_{[0,\alpha]} \eta(ds)$. It is
easy to check that $\pi_\alpha\pi_\beta\eta= \pi_\alpha\eta$ for
$\beta\ge\alpha\ge0$. Then the sequence of probability measures
$\{\mbf{Q}_{p,q}(\pi_k \eta,\cdot)\dvtx  k = 1, 2, \ldots\}$ induce a
consistent family of finite-dimensional distributions on the product
space $M_\infty:= \prod_{k=1}^\infty M[0,k]$. Let $\mbf{Q}$ be the
unique probability measure on $M_\infty$ determined by the family. Then
under $\mbf{Q}$ the canonical sequence $(X_1,X_2,\ldots)$ of $M_\infty$
converges almost surely to a random Radon measure $X$ on $[0,\infty)$,
which has distribution $\mbf{Q}_{p,q}(\eta,\cdot)$ on $\mcr
{M}[0,\infty)$ given by (\ref{4.18}). It is easy to show that
$\mbf{Q}_{p,q}(\eta,d\nu)$ is a probability kernel on the space
$\mcr{M}[0,\infty)$. The semigroup property of $\{\mbf{Q}_{p,q}\dvtx  q\ge
p\in T\}$ follows from (\ref{4.9}) and (\ref{4.10}).
\end{pf}

Since the state space $\mcr{M}[0,\infty)$ contains infinite measures, the
transition semigroup $\{\mbf{Q}_{p,q}\dvtx  q\ge p\in T\}$ defined by
(\ref{4.18}) does not fit exactly into the setup of the second section.
However, if $\{Z_q\dvtx  q\in T\}$ is a Markov process in $\mcr{M}[0,\infty)$
with transition semigroup $\{\mbf{Q}_{p,q}\dvtx  q\ge p\in T\}$, for each
$\alpha\ge0$, the restriction of $\{Z_q\dvtx  q\in T\}$ to $[0,\alpha]$
is an
inhomogeneous immigration superprocess with transition semigroup
$\{\mbf{Q}_{p,q}^\alpha\dvtx  q\ge p\in T\}$. Then we can think of the
original process $\{Z_q\dvtx  q\in T\}$ as an inhomogeneous immigration
superprocess with the \textit{extended state space} $\mcr{M}[0,\infty)$.
The model can be described intuitively as follows. The offspring born by
a ``particle'' at site $s\ge0$ at time $\theta\in T$ are spread over the
interval $[s,\infty)$ according to the law determined by
$\psi_\theta(s,\cdot)$. Thus the superprocess only involves a nonlocal
branching structure. The immigration rate is given by $\mu(d\theta)$ and
the immigrants coming at time $\theta\in T$ are distributed in
$[0,\infty)$ according to the law given by $u_\theta(0,\cdot)$. The
spatial motion of the immigration superprocess is trivial.

Suppose that $\{(X_t(q))_{t\ge0}\dvtx  q\in T\}$ is a Markov process with
transition semigroup $\{\mbf{P}_{p,q}\dvtx  q\ge p\in T\}$ defined by
(\ref{4.6}). We can identify the random path $(X_t(q))_{t\ge0}$ with the
absolutely continuous random measure $Z_q$ on $[0,\infty)$ with
$(X_t(q))_{t\ge0}$ as a density. By Theorem \ref{t4.5}, the
measure-valued process $\{Z_q\dvtx  q\in T\}$ is an immigration superprocess
with transition semigroup $\{\mbf{Q}_{p,q}\dvtx  q\ge p\in T\}$ defined by
(\ref{4.18}). Therefore we can naturally call $\{(X_t(q))_{t\ge0}\dvtx
q\in
T\}$ a \textit{path-valued branching process with immigration}. By
(\ref{4.5}) we have $X_q\ge X_p$ almost surely for $q\ge p\in T$. If
$\mu(q) = \mu$ independent of $q\in T$, we simply call $\{X_q\dvtx  q\ge
0\}$
a \textit{path-valued branching process}.

By (\ref{4.8}) or (\ref{4.12}) we have $u_{p,q}(s,f)\ge f(s)$ for any
$s\ge0$ and $f\in B^+[0,\infty)$ with compact support. Then (\ref{4.18})
implies that the set of infinite measures on $[0,\infty)$ is absorbing
for $\{\mbf{Q}_{p,q}\dvtx  q\ge p\in T\}$. Let $\{\mbf{Q}_{p,q}^\infty\dvtx
q\ge
p\in T\}$ denote the sub-Markov restriction of $\{\mbf{Q}_{p,q}\dvtx  q\ge
p\in T\}$ to the space $M[0,\infty)$ of finite measures on $[0,\infty)$.
If $f\in B^+[0,\infty)$ is bounded away from zero, we define
\[
u_q^\infty(s,f) = \lim_{\alpha\to\infty}u_q(s,f1_{[0,\alpha]}),\qquad
u_{p,q}^\infty(s,f) = \lim_{\alpha\to\infty} u_{p,q}(s,f1_{[0,\alpha]}).
\]
For an arbitrary $f\in B^+[0,\infty)$, define
\[
u_q^\infty(s,f) = \lim_{n\to\infty}u_q^\infty(s,f+1/n),\qquad
u_{p,q}^\infty(s,f) = \lim_{n\to\infty}u_{p,q}^\infty(s,f+1/n).
\]
By (\ref{4.15}) one can see $u_{p,q}^\infty(s):= u_{p,q}^\infty(s,f)$
solves
%
%e4.19 #&#
\begin{equation}
\label{4.19} u_{p,q}^\infty(s) = f(s) + \int
_p^q \psi_\theta^\infty
\bigl(s,u_{\theta,q}^\infty\bigr) \,d\theta,\qquad s\ge0, q\ge p\in T,
\end{equation}
where $\psi_\theta^\infty(s,f) = \zeta_\theta(u_\theta^\infty(s,f))$.
From (\ref{4.8}) we obtain
%
%e4.20 #&#
\begin{equation}
\label{4.20} u_{p,q}^\infty(s,f) = f(s) + \phi_{p,q}
\bigl(u_q^\infty(s,f)\bigr),\qquad s\ge0.
\end{equation}
It is easy to show that, for $f\in B^+[0,\infty)$,
%
%e4.21 #&#
\begin{equation}
\label{4.21} \int_{M[0,\infty)} e^{-\langle\nu,f\rangle} \mbf
{Q}_{p,q}^\infty(\eta,d\nu) = \exp\bigl\{-
\mu(p,q]u_q^\infty(0,f) - \bigl\langle\eta,u_{p,q}^\infty(
\cdot,f)\bigr\rangle\bigr\}.\hspace*{-30pt}
\end{equation}
To avoid the triviality of $\{\mbf{Q}_{p,q}^\infty\dvtx  q\ge p\in T\}$, we
need to assume $\phi_q(\lambda)\to\infty$ as $\lambda\to\infty$ for
every $q\in T$. In this case, we can define the right inverse
$\phi_q^{-1}$ of $\phi_q$ as in (\ref{2.20}). By (\ref{2.17}),
(\ref{2.19}) and (\ref{4.20}), we have
%
%e4.22 #&#
\begin{equation}
\label{4.22} u_q^\infty(s,\lambda) = \phi_q^{-1}(
\lambda),\qquad u_{p,q}^\infty(s,\lambda) = \phi_p
\bigl(\phi_q^{-1}(\lambda)\bigr),\qquad s\ge0, \lambda\ge0.\hspace*{-30pt}
\end{equation}

%th4.6 #&#
\begin{theorem}\label{t4.6} Suppose that $\phi_q(\lambda)\to\infty
$ as
$\lambda\to\infty$ for every $q\in T$. Let $S\subset T$ be an interval
not containing critical points of $\{\phi_q\dvtx  q\in T\}$. Then for any
$q\in S$ and $f\in B^+[0,\infty)$ there is a unique locally bounded
positive solution $(p,s)\mapsto u_{p,q}^\infty(s):= u_{p,q}^\infty(s,f)$
to (\ref{4.19}) on $S\times[0,\infty)$. Moreover,\vspace*{1pt} the sub-Markov
transition semigroup $\{\mbf{Q}_{p,q}^\infty\dvtx  q\ge p\in S\}$ defines an
inhomogeneous immigration superprocess in $M[0,\infty)$ with branching
mechanisms $\{-\psi_\theta^\infty\dvtx  \theta\in S\}$, immigration mechanisms
$\{u_\theta^\infty(0,\cdot)\dvtx  \theta\in S\}$ and immigration measure
$\mu$.
\end{theorem}
\begin{pf}
For any $s\ge0$ and $\theta\in S$, one can see by (\ref
{2.17}) that
the functional $f\mapsto u_\theta^\infty(s,f)$ on $B^+[0,\infty)$ is the
exponent of an infinitely divisible sub-probability measure carried by
$M[s,\infty)$. Then we have the representation
\[
u_\theta^\infty(s,f) = a_\theta^\infty(s) +
\bigl\langle\eta_\theta^\infty(s),f\bigr\rangle+ \int
_{M[s,\infty)^\circ} \bigl(1-e^{-\langle\nu,f\rangle}\bigr) H_\theta
^\infty(s,d
\nu),
\]
where $a_\theta^\infty(s)\ge0$, $\eta_\theta^\infty(s)\in
M[s,\infty)$
and $(1\land\langle\nu,1\rangle) H_\theta^\infty(s,d\nu)$ is a
finite measure on
$M[s,\infty)^\circ$. By the first equality in (\ref{4.22}) we get
$a_\theta^\infty(s) = u_\theta^\infty(s,0) = \phi_\theta^{-1}(0)$. It
follows that
\[
\bigl\langle\eta_\theta^\infty(s),1\bigr\rangle+ \int
_{M[s,\infty)^\circ
} \langle\nu,1\rangle H_\theta^\infty(s,d
\nu) = 1/\phi_\theta'\bigl(\phi_\theta^{-1}(0)
\bigr).
\]
The right-hand side is bounded on each compact subinterval of $S$. By
Theorem~\ref{t2.2}, the composed functional $f\mapsto\psi_\theta^\infty
(s,f) = \zeta_\theta(u_\theta^\infty(s,f))$ is the exponent of an
infinitely divisible sub-probability measure carried by $M[s,\infty)$.
Then $f\mapsto\psi_\theta^\infty(s,0) - \psi_\theta^\infty(s,f)$
can be
represented by a special form of (\ref{2.3}). That shows (\ref{4.19}) is
a special case of (\ref{2.7}). The desired result now follows in view of
(\ref{4.21}) and (\ref{4.17}) with $\alpha= \infty$.
\end{pf}

If $\phi_q(\lambda)\to\infty$ as $\lambda\to\infty$ for every
$q\in T$,
we can restrict $\{\mbf{P}_{p,q}\dvtx  q\ge p\in T\}$ to the space
$D^+_{\mrm{in}}[0,\infty)$ of integrable paths in $D^+[0,\infty)$ to get
a sub-Markov transition semigroup $\{\mbf{P}_{p,q}^\infty\dvtx  q\ge p\in
T\}$. This semigroup can also be regarded as a restriction of
$\{\mbf{Q}_{p,q}^\infty\dvtx  q\ge p\in T\}$. For $f\in B^+[0,\infty)$,
we have
%
%e4.23 #&#
\begin{eqnarray}
\label{4.23} &&\int_{D^+_{\mrm{in}}[0,\infty)} \exp\biggl\{-\int
_0^\infty f(s)w(s)\,ds \biggr\} \mbf{P}_{p,q}^\infty(
\eta,dw)
\nonumber\\[-8pt]\\[-8pt]
&&\qquad = \exp\biggl\{-\mu(p,q]u_q^\infty(0,f) - \int
_0^\infty u_{p,q}^\infty(s,f)
\eta(s) \,ds \biggr\}.
\nonumber
\end{eqnarray}

For an inhomogeneous immigration superprocess $\{Z_q\dvtx  q\in T\}$ with
transition semigroup $\{\mbf{Q}_{p,q}\dvtx  q\ge p\in T\}$ or
$\{\mbf{Q}_{p,q}^\infty\dvtx  q\ge p\in T\}$, we define its \textit{total mass
process} $\{\sigma(q)\dvtx  q\in T\}$ by $\sigma(q) = Z_q[0,\infty)$. For a
path-valued branching process with immigration $\{(X_t(q))_{t\ge0}\dvtx
q\in
T\}$ with transition semigroup $\{\mbf{P}_{p,q}\dvtx  q\ge p\in T\}$ or
$\{\mbf{P}_{p,q}^\infty\dvtx  q\ge p\in T\}$, its \textit{total mass process}
is defined as
\[
\sigma(q) = \int_0^\infty X_s(q)\,ds,\qquad
q\in T.
\]
We here think of $\{\sigma(q)\dvtx  q\in T\}$ as a process with state space
$\mbb{R}_+$ and cemetery $\infty$. In view of (\ref{4.21}), (\ref{4.22})
and (\ref{4.23}), we have
%
%th4.7 #&#
\begin{theorem}\label{t4.7} Suppose that $\phi_q(\lambda)\to\infty
$ as
$\lambda\to\infty$ for every $q\in T$. Then $\{\sigma(q)\dvtx  q\in T\}$ is
an inhomogeneous Markov process with transition semigroup $\{R_{p,q}\dvtx
q\ge p\in T\}$ such that, for $\lambda\ge0$,
%
%e4.24 #&#
\begin{equation}
\label{4.24} \int_{\mbb{R}_+} e^{-\lambda y} R_{p,q}(x,dy)
= \exp\bigl\{-x\phi_p\bigl(\phi^{-1}_q(
\lambda)\bigr) - \mu(p,q] \phi^{-1}_q(\lambda) \bigr\}.
\end{equation}
\end{theorem}

Before concluding this section, let us consider the admissible family of
branching mechanisms $\{\phi_q\dvtx  q\in\mbb{R}\}$ defined by
$\phi_q(\lambda) = \lambda^2 - 2q\lambda$ for $\lambda\ge0$. In this
special case, zero is the only critical point of the family $\{\phi_q\dvtx
q\in\mbb{R}\}$. Let $\{(X_t(q))_{t\ge0}\dvtx  q\in\mbb{R}\}$ be a
corresponding path-valued branching process. Let $\{\sigma(q)\dvtx  q\in
\mbb{R}\}$ be the process of total mass. By Theorem \ref{t4.7} one can
see that $\{\sigma(q)\dvtx  q\in\mbb{R}\}$ is an inhomogeneous Markov process
with transition semigroup $\{R_{p,q}\dvtx  q\ge p\in\mbb{R}\}$ defined by
%
%e4.25 #&#
\begin{equation}
\label{4.25} \int_{\mbb{R}_+} e^{-\lambda y} R_{p,q}(x,dy)
= \exp\bigl\{-xv_{p,q}(\lambda) \bigr\},\qquad \lambda\ge0,
\end{equation}
where
\[
v_{p,q}(\lambda) = \lambda+2(q-p) \bigl(\sqrt{q^2+
\lambda}+q\bigr).
\]
This process can be obtained from two homogeneous CB-processes by simple
transformations. For $t,\lambda\ge0$ let
\[
u_t^-(\lambda) = e^{-2t}\lambda+2e^{-t}
\bigl(1-e^{-t}\bigr) (\sqrt{1+\lambda}-1).
\]
It is easy to check that
\[
u_{t-s}^-(\lambda) = e^{2s}v_{-e^{-s},-e^{-t}}
\bigl(e^{-2t}\lambda\bigr),\qquad \lambda\ge0, t\ge s\in\mbb{R}.
\]
From this and (\ref{4.25}) one can see that $\{e^{-2t}\sigma(-e^{-t})\dvtx
t\in\mbb{R}\}$ is a homogeneous Markov process with transition semigroup
$(R^-_t)_{t\ge0}$ defined by
\[
\int_{\mbb{R}_+} e^{-\lambda y} R^-_t(x,dy) =
e^{-xu_t^-(\lambda)},\qquad \lambda\ge0.
\]
Moreover, we have
\[
\frac{d}{dt}u_t^-(\lambda) = - \phi_-\bigl(u_t^-(
\lambda)\bigr),
\]
where
\[
\phi_-(z) = 2z-2(\sqrt{1+z}-1).
\]
Then $\{e^{-2t}\sigma(-e^{-t})\dvtx  t\in\mbb{R}\}$ is actually a conservative
homogeneous CB-process in $[0,\infty)$ with branching mechanism $\phi_-$.
Similarly, one sees $\{e^{2t}\sigma(e^t)\dvtx\break  t\in\mbb{R}\}$ is a homogeneous
Markov process with transition semigroup $(R^+_t)_{t\ge0}$ defined by
\[
\int_{\mbb{R}_+} e^{-\lambda y} R^+_t(x,dy) =
e^{-xu_t^+(\lambda)},\qquad \lambda\ge0,
\]
where
\[
u_t^+(\lambda) = e^{2t}\lambda+2e^t
\bigl(e^t-1\bigr) (\sqrt{1+\lambda}+1).
\]
One can easily see that
\[
\frac{d}{dt}u_t^+(\lambda) = - \phi_+\bigl(u_t^+(
\lambda)\bigr),
\]
where
\[
\phi_+(z) = - 2z - 2(\sqrt{1+z}+1).
\]
Then $\{e^{2t}\sigma(e^t)\dvtx  t\in\mbb{R}\}$ is a CB-process with branching
mechanism $\phi_+$.

%%%%%%%%%%%% (Section 5) %%%%%%%%%%%%%%

%s5 #&#
\section{Construction by stochastic equations}\label{sec5}

In this section, we give a construction of the path-valued process
$\{(X_t(q))_{t\ge0}\dvtx  q\in T\}$ with transition semigroup
$\{\mbf{P}_{p,q}\dvtx  q\ge p\in T\}$ defined by (\ref{4.6}) as the solution
flow of a system of stochastic equations driven by time--space noises. We
shall assume $T = [0,\infty)$ or $[0,a]$ or $[0,a)$ for some $a>0$. This
specification of the index set is clearly not essential for the
applications. Let $\mu\in F(T)$, and let $\{\phi_q\dvtx  q\in T\}$ be an
admissible family of branching mechanisms, where $\phi_q$ is given by
(\ref{1.1}) with the parameters $(b,m) = (b_q,m_q)$ depending on $q\in
T$. Let $\mu(p,q] = \mu(q) - \mu(p)$ for $q\ge p\in T$, and let $m(dy,dz)$
be the measure on $T\times(0,\infty)$ defined by
%
%e5.1 #&#
\begin{equation}
\label{5.1} m\bigl([0,q]\times[c,d]\bigr) = m_q[c,d],\qquad q\in T, d>c>0.
\end{equation}
Let $\rho= \rho(s)$ be a locally bounded positive Borel function on
$[0,\infty)$, and let $\psi$ be an immigration mechanism given by
(\ref{2.21}).

Suppose that $(\iittOmega, \mcr{F}, \mcr{F}_t, \mbf{P})$ is a filtered
probability space satisfying the usual hypotheses. Let $W(ds,du)$ be an
$(\mcr{F}_t)$-white noise on $(0,\infty)^2$ based on the Lebesgue measure,
let $\tilde{N}_0(ds,dy,dz,du)$ be a compensated $(\mcr{F}_t)$-Poisson
random measure on $(0,\infty)\times T\times(0,\infty)^2$ with intensity
$dsm(dy,dz)\,du$ and let $N_1(ds,dz,du)$ be an $(\mcr{F}_t)$-Poisson random
measure on $(0,\infty)^3$ with intensity $dsn(dz)\,du$. Suppose that
$W(ds,du)$, $\tilde{N}_0(ds,dy,dz,du)$ and $N_1(ds,dz,du)$ are independent
of each other. For $q\in T$ it is easy to see that
\[
\tilde{N}(ds,dz,du):= \int_{\{0\le y\le q\}} \tilde{N}_0(ds,dy,dz,du)
\]
is a compensated Poisson random measure with intensity $dsm_q(dz)\,du$. By
Theorem \ref{t3.3} for every $q\in T$ there is a pathwise unique solution
to the stochastic equation
%
%e5.2 #&#
\begin{eqnarray}
\label{5.2} X_t(q) &=& \mu(q) - b_q\int
_0^t X_{s-}(q) \,ds + \sigma\int
_0^t\int_0^{X_{s-}(q)}
W(ds,du)
\nonumber
\\
&&{} + \int_0^t\int_{[0,q]}
\int_0^\infty\int_0^{X_{s-}(q)}
z \tilde{N}_0(ds,dy,dz,du)
\\
&&{} + h\int_0^t \rho(s) \,ds + \int
_0^t\int_0^\infty
\int_0^{\rho(s)} z N_1(ds,dz,du).
\nonumber
\end{eqnarray}
By Theorem \ref{t3.4} the solution $\{X_t(q)\dvtx  t\ge0\}$ is a CBI-process
with branching mechanism $\phi_q$, immigration mechanism $\psi$ and
immigration rate $\rho$.
%
%th5.1 #&#
\begin{theorem}\label{t5.1} The process $\{(X_t(q))_{t\ge0}\dvtx  q\in T\}
$ is a
path-valued branching process with immigration in $D^+[0,\infty)$ having
transition semigroup $\{\mbf{P}_{p,q}\dvtx\break   q\ge p\in T\}$ defined by
(\ref{4.6}).
\end{theorem}
\begin{pf}
We can rewrite equation (\ref{5.2}) into
%
%e5.3 #&#
\begin{eqnarray}
\label{5.3} X_t(q) &=& \mu(q) - h_q\int
_0^t X_{s-}(q) \,ds + \sigma\int
_0^t\int_0^{X_{s-}(q)}
W(ds,du)
\nonumber
\\
&&{} + h\int_0^t \rho(s) \,ds + \int
_0^t\int_0^q
\int_0^\infty\int_0^{X_{s-}(q)}
z N_0(ds,dy,dz,du)
\nonumber\\[-8pt]\\[-8pt]
&&{} + \int_0^t\int_{\{0\}}
\int_0^\infty\int_0^{X_{s-}(q)}
z \tilde{N}_0(ds,dy,dz,du)
\nonumber
\\
&&{} + \int_0^t\int_0^\infty
\int_0^{\rho(s)} z N_1(ds,dz,du),
\nonumber
\end{eqnarray}
where
\[
q\mapsto h_q:= b_q + \int_0^qd
\theta\int_0^\infty z n_\theta(dz) =
b_0 - \int_0^q
\beta_\theta \,d\theta
\]
is a decreasing function. Then, for $q\ge p\in T$, one can see by a simple
modification of Theorem 2.2 in \citet{DawLi12} that $X_t(q)\ge
X_t(p)$ for every $t\ge0$ with probability one. Let $\xi_t(p,q) = X_t(q)
- X_t(p)$ for $t\ge0$. From (\ref{5.3}) we have
%
%e5.4 #&#
\begin{eqnarray}
\label{5.4} \xi_t(p,q) &=& \mu(p,q] - b_q\int
_0^t \xi_{s-}(p,q) \,ds + \int
_p^q \beta_\theta \,d\theta\int
_0^t X_{s-}(p) \,ds
\nonumber
\\
&&{} + \sigma\int_0^t\int
_0^{\xi_{s-}(p,q)} W\bigl(ds,X_{s-}(p)+du\bigr)
\nonumber\\[-8pt]\\[-8pt]
&&{} + \int_0^t\int_{[0,q]}
\int_0^\infty\int_0^{\xi_{s-}(p,q)}
z \tilde{N}_0\bigl(ds,dy,dz,X_{s-}(p)+du\bigr)
\nonumber
\\
&&{} + \int_0^t \int_p^q
\int_0^\infty\int_0^{X_{s-}(p)}
z N_0(ds,dy,dz,du).
\nonumber
\end{eqnarray}
Here $W(ds,X_{s-}(p)+du)$ is a white noise based on the Lebesgue measure.
Note also that
\[
\int_{\{0\le y\le q\}} N_0\bigl(ds,dy,dz,X_{s-}(p)+du
\bigr)
\]
is a Poisson random measure with intensity $dsm_q(dz)\,du$, and
\[
\int_{\{p<y\le q\}} N_0(ds,dy,dz,du)
\]
is a Poisson random measure with intensity
\[
\int_{\{p<\theta\le q\}} \,dsn_\theta(dz)\,du \,d\theta.
\]
Clearly, the white noise and the two random measures are independent. By
Theorem \ref{t3.4}, conditioned upon $\{X_t(p)\dvtx  t\ge0\}$ the process
$\{\xi_t(p,q)\dvtx  t\ge0\}$ is a CBI-process with branching mechanism
$\phi_q$, immigration mechanism $\phi_{p,q}$ and immigration rate
$\{X_{t-}(p)\dvtx  t\ge0\}$. Conditioned upon $\{X_t(p)\dvtx  t\ge0\}$, the
process $\{\xi_t(p,q)\dvtx  t\ge0\}$ is clearly independent of the
$\sigma$-algebra generated by $\{X_t(v)\dvtx  t\ge0, v\in[0,p]\}$. Then
$\{(X_t(q))_{t\ge0}\dvtx  q\in T\}$ is a path-valued Markov process with
transition semigroup $\{\mbf{P}_{p,q}\dvtx  q\ge p\in T\}$.
\end{pf}
%
%th5.2 #&#
\begin{theorem}\label{t5.2} There is a positive function $(t,u)\mapsto C(t,u)$
on $[0,\infty)\times T$ bounded on compact sets so that, for any $t\ge0$
and $p\le q\le u\in T$,
%
%e5.5 #&#
\begin{eqnarray}
\label{5.5}
&&
\mbf{P} \Bigl\{\sup_{0\le s\le t}\bigl[X_s(q)-X_s(p)
\bigr] \Bigr\}\nonumber\\[-8pt]\\[-8pt]
&&\qquad \le C(t,u) \bigl\{\mu(p,q] + b_p-b_q
+ \sqrt{\mu(p,q]} + \sqrt{b_p-b_q} \bigr\}.
\nonumber
\end{eqnarray}
\end{theorem}
\begin{pf}
Since $\{X_t(p)\dvtx  t\ge0\}$ is a CBI-process, we see from
(\ref{3.2}) that $t\mapsto\mbf{P}[X_t(p)]$ is locally bounded. Let
$\{\xi_t(p,q)\dvtx  t\ge0\}$ be defined as in the last proof. By (\ref{5.4})
we have
\[
\mbf{P}\bigl[\xi_t(p,q)\bigr] = \mu(p,q] - b_q\int
_0^t \mbf{P}\bigl[\xi_s(p,q)\bigr]
\,ds + (b_p-b_q)\int_0^t
\mbf{P}\bigl[X_s(p)\bigr] \,ds.
\]
By Gronwall's inequality one can find a positive function $(t,u)\mapsto
C_0(t,u)$ on $[0,\infty)\times T$ bounded on compact sets so that, for
any $t\ge0$ and $p\le q\le u\in T$,
%
%e5.6 #&#
\begin{equation}
\label{5.6} \mbf{P}\bigl[\xi_t(p,q)\bigr] \le C_0(t,u)
\bigl\{\mu(p,q] + b_p-b_q \bigr\}.
\end{equation}
Applying Doob's inequality to the martingales in (\ref{5.4}), we obtain
\begin{eqnarray*}
\mbf{P} \Bigl\{\sup_{0\le s\le t}\xi_s(p,q) \Bigr\} &\le&
\mu(p,q] + 2\sigma\biggl(\int_0^t\mbf{P}\bigl[
\xi_s(p,q)\bigr]\,ds \biggr)^{1/2}
\\
&&{} + |b_q|\int_0^t \mbf{P}
\bigl[\xi_s(p,q)\bigr] \,ds + (b_p-b_q)\int
_0^t \mbf{P}\bigl[X_s(p)\bigr] \,ds
\\
&&{} + \int_1^\infty z m_q(dz)\int
_0^t \mbf{P}\bigl[\xi_s(p,q)\bigr]
\,ds
\\
&&{} + 2 \biggl(\int_0^t \mbf{P}\bigl[
\xi_s(p,q)\bigr]\,ds \int_0^1
z^2 m_q(dz) \biggr)^{1/2}.
\end{eqnarray*}
Then the desired estimate follows from (\ref{5.6}).
\end{pf}

Now let us consider a special admissible family of branching mechanisms.
Suppose that $\phi$ is a critical or supercritical branching mechanism
given by (\ref{1.1}) with $b\le0$. Let $T = T(\phi)$ be the set of
$q\ge
0$ so that
\[
\int_1^\infty ze^{qz}m(dz)< \infty.
\]
Then $T=[0,a]$ or $[0,a)$, where $a = \sup(T)$. We can define an
admissible family of branching mechanisms $\{\phi_q\dvtx  q\in T\}$ by
%
%e5.7 #&#
\begin{equation}
\label{5.7} \phi_q(\lambda) = \phi(\lambda-q) - \phi(-q),\qquad \lambda
\ge0,
\end{equation}
where the two terms on the right-hand side are defined using formula
(\ref{1.1}). Let $\{X_t(q)\dvtx  t\ge0,q\in T\}$ be the solution flow of
stochastic equation system (\ref{1.3}). By Theorem \ref{t5.1} we see that
$\{(X_t(q))_{t\ge0}\dvtx  q\in T\}$ is an inhomogeneous path-valued branching
process with transition semigroup $\{\mbf{P}_{p,q}\dvtx  q\ge p\in T\}$ given
by
%
%e5.8 #&#
\begin{equation}
\label{5.8}\qquad \int_{D^+[0,\infty)} e^{-\int_0^\infty f(s)w(s) \,ds} \mbf{P}_{p,q}(
\eta,dw) = \exp\biggl\{-\int_0^\infty
u_{p,q}(s,f)\eta(s) \,ds \biggr\},
\end{equation}
where $f\in B^+[0,\infty)$ has compact support, and $u_{p,q}(s,f)$ is given
by (\ref{4.8}). If $\phi(\lambda)\to\infty$ as $\lambda\to\infty
$, by
Theorem \ref{t4.7} the corresponding total mass process $\{\sigma(q)\dvtx
q\in
T\}$ is an inhomogeneous CB-process with transition semigroup $\{R_{p,q}\dvtx
q\ge p\in T\}$ given by
%
%e5.9 #&#
\begin{equation}
\label{5.9} \int_{\mbb{R}_+} e^{-\lambda y} R_{p,q}(x,dy)
= \exp\bigl\{-x\phi_p\bigl(\phi^{-1}_q(
\lambda)\bigr) \bigr\},\qquad \lambda\ge0.
\end{equation}
By Theorem \ref{t2.6} we have
%
%e5.10 #&#
\begin{equation}
\label{5.10} \mbf{P} \bigl[e^{-\lambda\sigma(q)}1_{\{\sigma(q)< \infty\}
}\bigr] =
e^{-\mu\phi^{-1}_q(\lambda)},\qquad \lambda\ge0, q\in T.
\end{equation}
It is simple to see that
\[
q\mapsto\phi^{-1}_q(0) = q + \phi^{-1}\bigl(
\phi(-q)\bigr)
\]
is continuous on $T$. Let $A = \inf\{q\in T\dvtx  \sigma(q) = \infty\}$
be the
explosion time of $\{\sigma(q)\dvtx  q\in T\}$. For any $q\in T$ we can let
$\lambda= 0$ in (\ref{5.10}) to obtain
%
%e5.11 #&#
\begin{equation}
\label{5.11} \mbf{P}\{A>q\} = \mbf{P}\bigl\{\sigma(q)<\infty\bigr\} =
e^{-\mu\phi^{-1}_q(0)}.
\end{equation}
This gives a characterization of the distribution of $A$.
%
%th5.3 #&#
\begin{theorem}\label{t5.3} Suppose that $\phi(\lambda)\to\infty$ as
$\lambda\to\infty$. Then for any $\theta\in T$, we have
%
%e5.12 #&#
\begin{equation}
\label{5.12} \mbf{P}\bigl[\sigma(\theta)1_{\{\sigma(\theta)<\infty\}
}\bigr] =
\frac{\mu e^{-\mu[\theta+ \phi^{-1}(\phi(-\theta))]}}{
\phi'(\phi^{-1}(\phi(-\theta)))}.
\end{equation}
\end{theorem}
\begin{pf}
Let $\lambda\ge0$ and $u=\phi_\theta^{-1}(\lambda)$. By
(\ref{5.10}) we have
\[
\mbf{P}\bigl[\sigma(\theta)e^{-\lambda\sigma(\theta)}1_{\{\sigma
(\theta)<\infty\}}\bigr] = -
\frac{d}{d\lambda}e^{-\mu\phi_\theta^{-1}(\lambda)} = \mu e^{-\mu\phi
_\theta^{-1}(\lambda)} \,\frac{d}{d\lambda}
\phi_\theta^{-1}(\lambda).
\]
From the relation $\phi_\theta(u) = \phi(u-\theta) - \phi(-\theta
)$, one
can see
\[
\phi_\theta^{-1}(\lambda) = \theta+ \phi^{-1}\bigl(
\lambda+\phi(-\theta)\bigr).
\]
It follows that
\[
\mbf{P}\bigl[\sigma(\theta)e^{-\lambda\sigma(\theta)}1_{\{\sigma
(\theta)<\infty\}}\bigr] =
\frac{\mu e^{-\mu[\theta+ \phi^{-1}(\lambda+\phi(-\theta))]}}{
\phi'(\phi^{-1}(\lambda+\phi(-\theta)))}.
\]
Then we get (\ref{5.12}) by letting $\lambda= 0$.
\end{pf}
%
%th5.4 #&#
\begin{theorem}\label{t5.4} Suppose that $\phi(\lambda)\to\infty$ as
$\lambda\to\infty$. Let $\theta\in[0,a)$, and let $G(\theta)$ be a
positive random variable measurable with respect to the $\sigma$-algebra
generated by $\{X_t(v)\dvtx  t\ge0, 0< v\le\theta\}$. Then we have
%
%e5.13 #&#
\begin{equation}
\label{5.13} \mbf{P}\bigl[G(\theta)|A=\theta\bigr] = \frac{\phi'(\phi
^{-1}(\phi(-\theta)))}{\mu e^{-\mu[\theta+
\phi^{-1}(\phi(-\theta))]}} \mbf{P}
\bigl[G(\theta)\sigma(\theta) 1_{\{\sigma(\theta)<\infty\}}\bigr].
\end{equation}
\end{theorem}
\begin{pf}
Since $q\mapsto\phi^{-1}_q(0)$ is continuous on $T$, for any
$q\in
(\theta,a)$ we can see by (\ref{5.9}) that
\[
\mbf{P}\bigl[G(\theta)1_{\{A>q\}}\bigr] = \mbf{P}\bigl[G(
\theta)1_{\{\sigma(q)<\infty\}}\bigr] = \mbf{P}\bigl[G(\theta)\exp\bigl
\{-\sigma(\theta)
\phi_\theta\bigl(\phi^{-1}_q(0)\bigr)\bigr\}\bigr].
\]
It is easy to see that
\[
\phi_\theta\bigl(\phi^{-1}_q(0)\bigr) =
\phi_\theta(\bar{q}+q) = \phi(\bar{q}+q-\theta) - \phi(-\theta),
\]
where $\bar{q} = \phi^{-1}(\phi(-q))$. By elementary calculations,
\[
\frac{d}{dq}\phi_\theta\bigl(\phi^{-1}_q(0)
\bigr) = \phi'(\bar{q}+q-\theta) \biggl(1-\frac{\phi'(-q)}{\phi'(\bar
{q})}
\biggr).
\]
It follows that
\[
- \frac{d}{dq}\mbf{P}\bigl[G(\theta)1_{\{A>q\}}\bigr]
\bigg|_{q=\theta+} = \bigl[\phi'(\bar{\theta})-\phi'(-
\theta)\bigr] \mbf{P}\bigl[G(\theta)\sigma(\theta) 1_{\{\sigma(\theta
)<\infty\}}\bigr],
\]
and hence
%
%e5.14 #&#
\begin{equation}
\label{5.14} \mbf{P}\bigl[G(\theta)|A=\theta\bigr] = \frac{\mbf
{P}[G(\theta)\sigma(\theta)1_{\{\sigma(\theta)<\infty\}
}]}{
\mbf{P}[\sigma(\theta) 1_{\{\sigma(\theta)<\infty\}}]}.
\end{equation}
Then we get (\ref{5.13}) from (\ref{5.12}) and (\ref{5.14}).
\end{pf}

%%%%%%%%%%%% (Section 6) %%%%%%%%%%%%%%

%s6 #&#
\section{A nonlocal branching superprocess}\label{sec6}

In this section, we consider a nonlocal branching superprocess defined
from the solution flow of (\ref{5.2}). We first assume $T = [0,a]$ for
some $a>0$. Let $\mu\in F(T)$, and let $\{\phi_q\dvtx  q\in T\}$ be an
admissible family of branching mechanisms, where $\phi_q$ is given by
(\ref{1.1}) with the parameters $(b,m) = (b_q,m_q)$ depending on $q\in T$.
Let $m(dy,dz)$ be the measure on $T\times(0,\infty)$ defined by
(\ref{5.1}). Let $\rho= \{\rho(t)\dvtx  t\ge0\}$ be a locally bounded
positive Borel function on $[0,\infty)$. Let $\psi$ be an immigration
mechanism given by (\ref{2.21}). Let $X(q) = \{X_t(q)\dvtx  t\ge0\}$ be the
solution of (\ref{5.2}) for $q\in T$. Then the path-valued Markov process
$\{X(q)\dvtx  q\in T\}$ has transition semigroup $\{\mbf{P}_{p,q}\dvtx  q\ge
p\in
T\}$ defined by (\ref{4.6}). Let $Q_T$ denote the set of rationals in $T$.
For any $t\ge0$ we define the random function $Y_t\in F(T)$ by $Y_t(a) =
X_t(a)$ and
%
%e6.1 #&#
\begin{equation}
\label{6.1} Y_t(q) = \inf\bigl\{X_t(u)\dvtx  u\in
Q_T\cap(q,a]\bigr\},\qquad 0\le q<a.
\end{equation}
Similarly, for any $t>0$, define $Z_t\in F(T)$ by $Z_t(a) = X_{t-}(a)$ and
%
%e6.2 #&#
\begin{equation}
\label{6.2} Z_t(q) = \inf\bigl\{X_{t-}(u)\dvtx  u\in
Q_T\cap(q,a]\bigr\},\qquad 0\le q<a.
\end{equation}
By Theorem \ref{t5.2}, for each $q\in T$ we have
%
%e6.3 #&#
\begin{equation}
\label{6.3} \mbf{P}\bigl\{Y_t(q)=X_t(q) \mbox{ and }
Z_t(q)=X_{t-}(q) \mbox{ for all } t\ge0\bigr\} = 1.
\end{equation}
Consequently, for every $q\in T$ the process $\{Y_t(q)\dvtx  t\ge0\}$ is
almost surely c\`adl\`ag and solves (\ref{5.2}), so it is a CBI-process
with branching mechanism $\phi_q$, immigration mechanism $\psi$ and
immigration rate $\rho$. In view of (\ref{4.3}) and (\ref{4.4}), for every
$q\in T$ we almost surely have
%
%e6.4 #&#
\begin{eqnarray}
\label{6.4} Y_t(q) &=& \mu(q) + A_t + \sigma\int
_0^t\int_0^{Y_{s-}(q)}
W(ds,du)
\nonumber
\\
&&{} - b_0 \int_0^t
Y_{s-}(q) \,ds + \int_0^q
\beta_\theta \,d\theta\int_0^t
Y_{s-}(q) \,ds
\nonumber\\[-8pt]\\[-8pt]
&&{} + \int_0^t\int_{\{0\}}
\int_0^\infty\int_0^{Y_{s-}(q)}
z \tilde{N}_0(ds,dy,dz,du)
\nonumber
\\
&&{} + \int_0^t\int_0^q
\int_0^\infty\int_0^{Y_{s-}(q)}
z N_0(ds,dy,dz,du),
\nonumber
\end{eqnarray}
where
\[
A_t = h\int_0^t \rho(s) \,ds + \int
_0^t\int_0^\infty
\int_0^{\rho
(s)} z N_1(ds,dz,du).
\]

For $t\ge0$ let $Y_t(dx)$ and $Z_t(dx)$ denote the random measures on $T$
induced by the random functions $Y_t$ and $Z_t\in F(T)$, respectively. For
any $f\in C^1(T)$ one can use Fubini's theorem to see
%
%e6.5 #&#
\begin{equation}
\label{6.5} \langle Y_t,f\rangle= f(a)Y_t(a)-\int
_0^a f'(q)Y_t(q)\,dq.
\end{equation}
Fix an integer $n\ge1$ and let $q_i = ia/2^n$ for $i=0,1,\ldots,2^n$. By
(\ref{6.3}) and (\ref{6.4}) it holds almost surely that
%
%e6.6 #&#
\begin{eqnarray}
\label{6.6}
&&\sum_{i=1}^{2^n}f'(q_i)Y_t(q_i)\nonumber\\
&&\qquad= \sum_{i=1}^{2^n}f'(q_i)
\mu(q_i) + \sigma\sum_{i=1}^{2^n}f'(q_i)
\int_0^t\int_0^{Z_s(q_i)}
W(ds,du)
\nonumber
\\
&&\qquad\quad{} + A_t\sum_{i=1}^{2^n}f'(q_i)
- b_0\sum_{i=1}^{2^n}f'(q_i)
\int_0^t Z_s(q_i) \,ds
\nonumber
\\
&&\qquad\quad{} + \sum_{i=1}^{2^n}f'(q_i)
\int_0^{q_i} \beta_\theta \,d\theta\int
_0^t ds \int_{[0,q_i]}Z_s(dx)
\nonumber
\\
&&\qquad\quad{} + \sum_{i=1}^{2^n}f'(q_i)
\int_0^t\int_{\{0\}}\int
_0^\infty\int_0^{Z_s(q_i)}
z \tilde{N}_0(ds,dy,dz,du)
\nonumber\\
&&\qquad\quad{} + \sum_{i=1}^{2^n}f'(q_i)
\int_0^t \int_0^{q_i}
\int_0^\infty\int_0^{Z_s(q_i)}
z N_0(ds,dy,dz,du)
\\
&&\qquad= \sum_{i=1}^{2^n}f'(q_i)
\mu(q_i) + \sigma\int_0^t\int
_0^{Z_s(a)} F_n(s,0,u) W(ds,du)
\nonumber
\\
&&\qquad\quad{} + A_t\sum_{i=1}^{2^n}f'(q_i)
- b_0\int_0^t \Biggl[\sum
_{i=1}^{2^n}f'(q_i)Z_s(q_i)
\Biggr] \,ds
\nonumber
\\
&&\qquad\quad{} + \int_0^t ds \int_T
Z_s(dx) \int_0^a
F_n(s,x\vee\theta,0)\beta_\theta \,d\theta
\nonumber
\\
&&\qquad\quad{} + \int_0^t\int_{\{0\}}
\int_0^\infty\int_0^{Z_s(a)}
zF_n(s,0,u) \tilde{N}_0(ds,dy,dz,du)
\nonumber
\\
&&\qquad\quad{} + \int_0^t \int_0^a
\int_0^\infty\int_0^{Z_s(a)}
zF_n(s,y,u) N_0(ds,dy,dz,du),
\nonumber
\end{eqnarray}
where
\[
F_n(s,y,u) = \sum_{i=1}^{2^n}f'(q_i)
1_{\{y\le q_i\}}1_{\{u\le
Z_s(q_i)\}}.
\]
By the right continuity of $q\mapsto Z_s(q)$ it is not hard to see that,
as $n\to\infty$,
%
%e6.7 #&#
\begin{equation}
\label{6.7} 2^{-n}F_n(s,y,u)\to F(s,y,u):= \int
_y^a 1_{\{u\le Z_s(q)\}} f'(q)\,dq.
\end{equation}
Then we can multiply (\ref{6.6}) by $2^{-n}$ and let $n\to\infty$ to see,
almost surely,
%
%e6.8 #&#
\begin{eqnarray}
\label{6.8}
&& \int_0^af'(q)Y_t(q)\,dq\nonumber\\
&&\qquad= \int_0^af'(q)\mu(q)\,dq +
\sigma\int_0^t\int_0^{Z_s(a)}
F(s,0,u) W(ds,du)
\nonumber
\\
&&\qquad\quad{} + A_t\int_0^a
f'(q) \,dq - b_0\int_0^t
ds\int_0^af'(q)Z_s(q)\,dq
\nonumber\\[-8pt]\\[-8pt]
&&\quad\qquad{} + \int_0^t ds\int_T
Z_s(dx) \int_0^a F(s,x\vee
\theta,0)\beta_\theta \,d\theta
\nonumber\\
&&\qquad\quad{} + \int_0^t\int_{\{0\}}
\int_0^\infty\int_0^{Z_s(a)}
zF(s,0,u) \tilde{N}_0(ds,dy,dz,du)
\nonumber
\\
&&\qquad\quad{} + \int_0^t \int_0^a
\int_0^\infty\int_0^{Z_s(a)}
zF(s,y,u) N_0(ds,dy,dz,du).
\nonumber
\end{eqnarray}
From (\ref{6.4}), (\ref{6.5}) and (\ref{6.8}) it follows that, almost
surely,
%
%e6.9 #&#
\begin{eqnarray}
\label{6.9} \langle Y_t,f\rangle&=& \langle\mu,f\rangle+
f(0)A_t + \sigma\int_0^t \int
_0^{Z_s(a)} \bigl[f(a)-F(s,0,u)\bigr] W(ds,du)
\nonumber
\\
&&{} - b_0\int_0^t \langle
Z_s,f\rangle \,ds + \int_0^t ds\int
_T Z_s(dx)\int_0^a
f(x\vee\theta) \beta_\theta \,d\theta
\nonumber\\[-8pt]\\[-8pt]
&&{} + \int_0^t\int_{\{0\}}
\int_0^\infty\int_0^{Z_s(a)}
z\bigl[f(a)-F(s,0,u)\bigr] \tilde{N}_0(ds,dy,dz,du)
\nonumber
\\
&&{} + \int_0^t\int_0^a
\int_0^\infty\int_0^{Z_s(a)}
z\bigl[f(a)-F(s,y,u)\bigr] N_0(ds,dy,dz,du).
\nonumber
\end{eqnarray}

%th6.1 #&#
\begin{theorem}\label{t6.1} The measure-valued process $\{Y_t\dvtx  t\ge0\}$
has a
c\`adl\`ag modification.
\end{theorem}
\begin{pf}
By (\ref{6.9}) one can see $\{\langle Y_t,f\rangle\dvtx  t\ge0\}$
has a c\`adl\`ag
modification for every $f\in C^1(T)$. Let $\mcr{U}$ be the countable set
of polynomials having rational coefficients. Then $\mcr{U}$ is uniformly
dense in both $C^1(T)$ and $C(T)$. For $f\in\mcr{U}$, let $\{Y_t^*(f)\dvtx
t\ge0\}$ be a c\`adl\`ag modification of $\{\langle Y_t,f\rangle\dvtx
t\ge0\}$. By
removing a null set from $\iittOmega$ if it is necessary, we obtain a
c\`adl\`ag process $\{Y_t^*\dvtx  t\ge0\}$ of rational linear functionals on
$\mcr{U}$, which can immediately be extended to a c\`adl\`ag process of
real linear functionals on $C(T)$. By Riesz's representation, the latter
determines a measure-valued process, which is clearly a c\`adl\`ag
modification of $\{Y_t\dvtx  t\ge0\}$.
\end{pf}
%
%th6.2 #&#
\begin{theorem}\label{t6.2} The c\`adl\`ag modification of $\{Y_t\dvtx
t\ge0\}
$ is
the unique solution of the following martingale problem: for every
$G\in
C^2(\mbb{R})$ and $f\in C(T)$,
%
%e6.10 #&#
\begin{eqnarray}
\label{6.10} G\bigl(\langle Y_t,f\rangle\bigr) &=& G\bigl(\langle
\mu,f\rangle\bigr) + \int_0^t G'
\bigl(\langle Y_s,f\rangle\bigr)\,ds \int_T
Y_s(dx) \int_T f(x\vee\theta)
\beta_\theta \,d\theta
\nonumber
\\
&&{} - b_0\int_0^t
G'\bigl(\langle Y_s,f\rangle\bigr)\langle
Y_s,f\rangle \,ds + \frac{1}{2} \sigma^2 \int
_0^t G''\bigl(
\langle Y_s,f\rangle\bigr)\bigl\langle Y_s,f^2
\bigr\rangle \,ds
\nonumber
\\
&&{} + \int_0^tds\int_T
Y_s(dx)\int_0^\infty\bigl[G\bigl(
\langle Y_s,f\rangle+ zf(x)\bigr)
\nonumber
\\
&&\hspace*{108pt}{} - G\bigl(\langle Y_s,f\rangle\bigr) - zf(x)G'
\bigl(\langle Y_s,f\rangle\bigr) \bigr]m_0(dz)
\nonumber\\
&&{} + \int_0^tds\int_T
Y_s(dx)\int_T \,d\theta\int
_0^\infty\bigl[G\bigl(\langle Y_s,f
\rangle+ zf(x\vee\theta)\bigr)
\\
&&\hspace*{184pt}{} - G\bigl(\langle Y_s,f\rangle\bigr) \bigr]n_\theta(dz)\nonumber\\
&&{}+ hf(0)\int_0^t G'\bigl(\langle
Y_s,f\rangle\bigr)\rho(s) \,ds
\nonumber
\\
&&{} + \int_0^t \rho(s)\,ds \int
_0^\infty\bigl[G\bigl(\langle Y_s,f
\rangle+ zf(0)\bigr) - G\bigl(\langle Y_s,f\rangle\bigr) \bigr] n(dz)
\nonumber
\\
&&{} + \mbox{local mart.}
\nonumber
\end{eqnarray}
\end{theorem}
\begin{pf}
We first assume $f\in C^1(T)$. By (\ref{6.9}) and It\^o's formula,
we get
\begin{eqnarray*}
\hspace*{-4pt}&&
G\bigl(\langle Y_t,f\rangle\bigr) \\
\hspace*{-4pt}&&\qquad= G\bigl(\langle\mu,f\rangle
\bigr) - b_0\int_0^t
G'\bigl(\langle Z_s,f\rangle\bigr)\langle
Z_s,f\rangle \,ds
\\
\hspace*{-4pt}&&\qquad\quad{} + \frac{1}{2}\sigma^2 \int_0^tds
\int_0^{Z_s(a)} G''\bigl(
\langle Z_s,f\rangle\bigr) \bigl[f(a) - F(s,0,u)\bigr]^2
\,du
\\
\hspace*{-4pt}&&\qquad\quad{} + \int_0^t G'\bigl(\langle
Z_s,f\rangle\bigr)\,ds\int_T
Z_s(dx)\int_T f(x\vee\theta)
\beta_\theta \,d\theta
\\
\hspace*{-4pt}&&\qquad\quad{} + \int_0^tds\int_0^{Z_s(a)}
\,du\int_0^\infty\bigl[G\bigl(\langle
Z_s,f\rangle+ z\bigl[f(a)-F(s,0,u)\bigr]\bigr)
\\
\hspace*{-4pt}&&\qquad\quad\hspace*{109pt}{}- G\bigl(\langle Z_s,f\rangle\bigr)\\
\hspace*{-4pt}&&\hspace*{109pt}\qquad\quad{}  - z\bigl[f(a)-F(s,0,u)
\bigr]G'\bigl(\langle Z_s,f\rangle\bigr)
\bigr]m_0(dz)
\\
\hspace*{-4pt}&&\qquad\quad{} + \int_0^tds\int_0^{Z_s(a)}
du\int_T \,d\theta\int_0^\infty
\bigl[G\bigl(\langle Z_s,f\rangle+ z\bigl[f(a)-F(s,\theta,u)\bigr]
\bigr)
\\
\hspace*{-4pt}&&\hspace*{230pt}\qquad\quad{} - G\bigl(\langle Z_s,f\rangle\bigr)
\bigr]n_\theta(dz)\\
\hspace*{-4pt}&&\qquad\quad{}+ hf(0)\int_0^t G'\bigl(\langle
Z_s,f\rangle\bigr)\rho(s) \,ds
\\
\hspace*{-4pt}&&\qquad\quad{} + \int_0^t \rho(s)\,ds \int
_0^\infty\bigl[G\bigl(\langle Z_s,f
\rangle+ zf(0)\bigr) - G\bigl(\langle Z_s,f\rangle\bigr) \bigr] n(dz)
\\
\hspace*{-4pt}&&\qquad\quad{} + \mbox{local mart.}
\end{eqnarray*}
For $s,u>0$ let $Z_s^{-1}(u) = \inf\{q\ge0\dvtx  Z_s(q)>u\}$. It is easy to
see that $\{q\ge0\dvtx  u\le Z_s(q)\} = [Z_s^{-1}(u),\infty)$, except for at
most countably many $u>0$. Then in the above we can replace $f(a) -
F(s,\theta,u)$ by
\[
f(a)-\int_\theta^a 1_{\{Z_s^{-1}(u)\le q\}}
f'(q)\,dq = f\bigl(Z_s^{-1}(u)\vee\theta
\bigr).
\]
It follows that
\begin{eqnarray*}
&&
G\bigl(\langle Y_t,f\rangle\bigr) \\
&&\qquad= G\bigl(\langle\mu,f\rangle
\bigr) - b_0\int_0^t
G'\bigl(\langle Z_s,f\rangle\bigr)\langle
Z_s,f\rangle \,ds
\\
&&\qquad\quad{} + \frac{1}{2}\sigma^2 \int_0^tds
\int_0^{Z_s(a)} G''\bigl(
\langle Z_s,f\rangle\bigr) f\bigl(Z_s^{-1}(u)
\bigr)^2 \,du
\\
&&\qquad\quad{} + \int_0^t G'\bigl(\langle
Z_s,f\rangle\bigr)\,ds\int_TZ_s(dx)
\int_T f(x\vee\theta) \beta_\theta \,d\theta
\\
&&\qquad\quad{} + \int_0^tds\int_0^{Z_s(a)}
du\int_0^\infty\bigl[G\bigl(\langle
Z_s,f\rangle+ zf\bigl(Z_s^{-1}(u)\bigr)\bigr)
\\
&&\qquad\quad\hspace*{107pt}{} - G\bigl(\langle Z_s,f\rangle\bigr)\\
&&\qquad\quad\hspace*{107pt}{} - zf\bigl(Z_s^{-1}(u)
\bigr)G'\bigl(\langle Z_s,f\rangle\bigr)
\bigr]m_0(dz)
\\
&&\qquad\quad{} + \int_0^tds\int_0^{Z_s(a)}
du\int_T \,d\theta\int_0^\infty
\bigl[G\bigl(\langle Z_s,f\rangle+ zf\bigl(Z_s^{-1}(u)
\vee\theta\bigr)\bigr)
\\
&&\qquad\quad\hspace*{209.5pt}{} - G\bigl(\langle Z_s,f\rangle\bigr)
\bigr]n_\theta(dz)\\
&&\qquad\quad{}+ hf(0)\int_0^t G'\bigl(\langle
Z_s,f\rangle\bigr)\rho(s) \,ds
\\
&&\qquad\quad{} + \int_0^t \rho(s)\,ds \int
_0^\infty\bigl[G\bigl(\langle Z_s,f
\rangle+ zf(0)\bigr) - G\bigl(\langle Z_s,f\rangle\bigr) \bigr] n(dz)
\\
&&\qquad\quad{} + \mbox{local mart.}
\\
&&\qquad= G\bigl(\langle\mu,f\rangle\bigr) + \int_0^t
G'\bigl(\langle Z_s,f\rangle\bigr)\,ds\int
_TZ_s(dx) \int_T f(x
\vee\theta) \beta_\theta \,d\theta
\\
&&\qquad\quad{} - b_0\int_0^t
G'\bigl(\langle Z_s,f\rangle\bigr)\langle
Z_s,f\rangle \,ds + \frac{1}{2}\sigma^2 \int
_0^t G''\bigl(
\langle Z_s,f\rangle\bigr)\bigl\langle Z_s,f^2
\bigr\rangle \,ds
\\
&&\qquad\quad{} + \int_0^tds\int_T
Z_s(dx)\int_0^\infty\bigl[G\bigl(
\langle Z_s,f\rangle+ zf(x)\bigr)
\\
&&\hspace*{109pt}\qquad\quad{} - G\bigl(\langle Z_s,f\rangle\bigr) - zf(x)G'
\bigl(\langle Z_s,f\rangle\bigr) \bigr]m_0(dz)
\\
&&\qquad\quad{} + \int_0^tds\int_T
Z_s(dx)\int_T \,d\theta\int
_0^\infty\bigl[G\bigl(\langle Z_s,f
\rangle+ zf(x\vee\theta)\bigr)
\\
&&\hspace*{186pt}\qquad\quad{} - G\bigl(\langle Z_s,f\rangle\bigr)
\bigr]n_\theta(dz)\\
&&\qquad\quad{}+ hf(0)\int_0^t G'\bigl(\langle
Z_s,f\rangle\bigr)\rho(s) \,ds
\\
&&\qquad\quad{} + \int_0^t \rho(s)\,ds \int
_0^\infty\bigl[G\bigl(\langle Z_s,f
\rangle+ zf(0)\bigr) - G\bigl(\langle Z_s,f\rangle\bigr) \bigr] n(dz)
\\
&&\qquad\quad{} + \mbox{local mart.}
\end{eqnarray*}
For each $q\in T$ the c\`adl\`ag process $\{X_t(q)\dvtx  t\ge0\}$ has at most
countably many discontinuity points $A_q:= \{t>0\dvtx  Y_{t-}(q)\neq
Y_t(q)\}$. In view of (\ref{6.1}) and (\ref{6.2}), we have $Z_t(q) =
Y_t(q)$ for all $q\in T$ and $t\in B:= A_a^c\cap(\bigcap_{u\in Q_T}A_u^c)$.
Here $B\subset[0,\infty)$ is a set with full Lebesgue measure. Then we
have (\ref{6.10}) for $f\in C^1(T)$. For an arbitrary $f\in C(T)$, we get
(\ref{6.10}) by an approximation argument. The uniqueness (in
distribution) of the solution to the martingale problem follows by a
modification of the proof of Theorem 7.13 in \citet{Li11}.
\end{pf}

The martingale problem (\ref{6.10}) is essentially a special case of the
one given in Theorem 10.18 of \citet{Li11}; see also Theorem 9.18 of \citet{Li11}. Let $f\mapsto\iittPsi(\cdot,f)$ be the operator on $C^+(T)$ defined
by
%
%e6.11 #&#
\begin{equation}
\label{6.11}\quad \iittPsi(x,f) = \int_T f(x\vee\theta)
\beta_\theta \,d\theta+ \int_T \,d\theta\int
_0^\infty\bigl(1-e^{-zf(x\vee
\theta)}
\bigr)n_\theta(dz).
\end{equation}
By modifying the proof of Theorem \ref{t3.4} one can show the following:
%
%th6.3 #&#
\begin{theorem}\label{t6.3} The solution $\{Y_t\dvtx  t\ge0\}$ of the martingale
problem (\ref{6.10}) is an immigration superprocess with transition
semigroup $(Q_t)_{t\ge0}$ defined by
%
%e6.12 #&#
\begin{equation}
\label{6.12}\qquad \int_{M(T)}e^{-\langle\nu,f\rangle}Q_t(\mu,d
\nu) = \exp\biggl\{-\langle\mu,V_tf\rangle- \int_0^t
\psi\bigl(V_sf(0)\bigr)\rho(s) \,ds \biggr\},
\end{equation}
where $f\in C^+(T)$, and $t\mapsto V_tf$ is the unique locally bounded
positive solution of
%
%e6.13 #&#
\begin{equation}
\label{6.13}\qquad V_tf(x) = f(x) - \int_0^t
\bigl[\phi_0\bigl(V_sf(x)\bigr) - \iittPsi(x,V_sf)
\bigr]\,ds,\qquad t\ge0, x\in T.
\end{equation}
\end{theorem}

The branching mechanism of the immigration superprocess $\{Y_t\dvtx  t\ge0\}$
has \textit{local part} $(x,f)\mapsto\phi_0(f(x))$ and \textit{nonlocal
part} $(x,f)\mapsto\iittPsi(x,f)$; see Example 2.5 in \citet{Li11}. The
process has immigration mechanism $f\mapsto\psi(f(0))$ and immigration
rate $\rho= \{\rho(s)\dvtx  s\ge0\}$. Then the immigrants only come at the
origin. The spatial motion in this model is trivial. Heuristically, when
an infinitesimal particle dies at site $x\in T$, some offspring are born
at this site according to the local branching mechanism and some are born
in the interval $(x,a]$ according to the nonlocal branching mechanism.
Therefore the branching of an infinitesimal particle located at $x\in T$
does not make any influence on the population in the interval $[0,x)$.
This explains the Markov property of the path-valued process
$\{(Y_t(q))_{t\ge0}\dvtx  q\in T\}$.

The cumulant semigroup $(V_t)_{t\ge0}$ can also be defined by a
differential evolution equation. In fact, by Theorem 7.11 of \citet{Li11}, for
any $f\in C^+(T)$, the integral equation (\ref{6.13}) is equivalent to
%
%e6.14 #&#
\begin{equation}
\label{6.14} \cases{ %
\displaystyle \frac{dV_tf}{dt}(x) = - \phi_0
\bigl(V_tf(x)\bigr) + \iittPsi(x,V_tf), &\quad $t\ge0, x\in T$,
\vspace*{2pt}\cr
V_0f(x) = f(x), &\quad $x\in T$.}
\end{equation}
Then the transition semigroup $(Q_t)_{t\ge0}$ can also be defined by
(\ref{6.12}) for $f\in C^+(T)$ with $t\mapsto V_tf$ being the unique
locally bounded positive solution of (\ref{6.14}).
%
%th6.4 #&#
\begin{theorem}\label{t6.4} Let $Y = (\iittOmega, \mcr{G}, \mcr{G}_t, Y_t,
\mbf{Q}_\mu)$ be any c\`adl\`ag immigration superprocess with transition
semigroup $(Q_t)_{t\ge0}$ defined by (\ref{6.12}) and (\ref{6.13}). Then
under $\mbf{Q}_\mu$ for every $q\in T$ the process $\{Y_t[0,q]\dvtx  t\ge
0\}$
has a c\`adl\`ag version, and $\{(Y_t[0,q])_{t\ge0}\dvtx  q\in T\}$ is a
path-valued branching process with immigration with transition semigroup
$\{\mbf{P}_{p,q}\dvtx  q\ge p\in T\}$ defined by (\ref{4.6}).
\end{theorem}
\begin{pf}
By Theorem \ref{t6.2}, one can see that for each $q\in T$ the
restriction of $\{Y_t\dvtx  t\ge0\}$ to $[0,q]$ is also an immigration
superprocess with state space $M[0,q]$. In particular, the process
$\{Y_t[0,q]\dvtx  t\ge0\}$ has a c\`adl\`ag version. Clearly, the
finite-dimensional distributions of the path-valued process
$\{(Y_t[0,q])_{t\ge0}\dvtx  q\in T\}$ are uniquely determined by the initial
state $\mu\in M(T)$ and transition semigroup $(Q_t)_{t\ge0}$. Then
$\{(Y_t[0,q])_{t\ge0}\dvtx  q\in T\}$ has identical finite-dimensional
distributions with the process $\{(Y_t(q))_{t\ge0}\dvtx  q\in T\}$ defined by
(\ref{6.4}). Since $\{(Y_t(q))_{t\ge0}\dvtx  q\in T\}$ is a Markov process
with transition semigroup $\{\mbf{P}_{p,q}\dvtx  q\ge p\in T\}$, so is
$\{(Y_t[0,q])_{t\ge0}\dvtx  q\in T\}$.
\end{pf}

If $T = [0,\infty)$ or $[0,a)$ for some $a>0$, we may apply the above
results to the interval $[0,q]\subset T$ for $q\in T$. Then for each
$q\in
T$, there is a immigration superprocess $\{Y_t^q\dvtx  t\ge0\}$ in $M[0,q]$.
Those processes determine a nonlocal branching immigration superprocess
$\{Y_t\dvtx  t\ge0\}$ in $\mcr{M}(T)$, the space of Radon measures on $T$
furnished with the topology of vague convergence. The results established
in this section hold for this process with obvious modifications.

%%%%%%%%%%%% (Section 7) %%%%%%%%%%%%%%

%s7 #&#
\section{The excursion law}\label{sec7}

In this section we assume $T=[0,a]$ for some $a>0$. However, the results
obtained here can be modified to the case $T = [0,a)$ or $[0,\infty)$,
obviously. Let $\{\phi_q\dvtx  q\in T\}$ be an admissible family of branching
mechanisms, where $\phi_q$ is given by (\ref{1.1}) with the parameters
$(b,m) = (b_q,m_q)$ depending on $q\in T$. In addition, we assume
$\phi_0'(\lambda)\to\infty$ as $\lambda\to\infty$. By
Theorem \ref{t6.3}, we can define the transition semigroup
$(Q_t)_{t\ge
0}$ of a nonlocal branching superprocess by
%
%e7.1 #&#
\begin{equation}
\label{7.1} \int_{M(T)}e^{-\langle\nu,f\rangle}Q_t(\mu,d
\nu) = \exp\bigl\{-\langle\mu,V_tf\rangle\bigr\},\qquad f\in C^+(T),
\end{equation}
where $t\mapsto V_tf$ is the unique locally bounded positive solution of
(\ref{6.13}). Let $(Q_t^\circ)_{t\ge0}$ denote the restriction of the
semigroup to $M(T)^\circ$.
%
%th7.1 #&#
\begin{theorem}\label{t7.1} The cumulant semigroup of $(V_t)_{t\ge0}$ in
(\ref{7.1}) admits the representation
%
%e7.2 #&#
\begin{equation}
\label{7.2} V_tf(x) = \int_{M(T)^\circ}
\bigl(1-e^{-\langle\nu,f\rangle}\bigr)L_t(x,d\nu),\qquad t>0, x\in T,
\end{equation}
where $(L_t(x,\cdot))_{t>0}$ is a $\sigma$-finite entrance law for
$(Q^\circ_t)_{t\ge0}$.
\end{theorem}
\begin{pf}
We need a modification of the characterization (\ref{6.14}) of the
cumulant semigroup. Let us consider a jump process $\xi$ in $T$ with
generator $A$ defined by
\[
Af(x) = \int_0^a \bigl(f(q)-f(x)\bigr)
\gamma(dq),\qquad x\in T, f\in C(T),
\]
where
\[
\gamma(dq) = \beta_q\,dq + \int_{\{0<z<\infty\}}
zn_q(dz)\,dq.
\]
Let $\phi_*(\lambda) = \gamma[0,a]\lambda+ \phi_0(\lambda)$, and let
$f\mapsto\iittPsi_*(\cdot,f)$ be the operator on $C^+(T)$ defined by
\[
\iittPsi_*(x,f) = \int_0^a\int
_0^\infty\bigl[e^{-zf(x\vee y)} - 1 + zf(x\vee y)
\bigr] m(dy,dz).
\]
Now the first equation in (\ref{6.14}) can be rewritten as
\[
\frac{dV_tf}{dt}(x) = AV_tf(x) - \phi_*\bigl(V_tf(x)
\bigr) - \iittPsi_*(x,V_tf).
\]
Then we may think of $(V_t)_{t\ge0}$ as the cumulant semigroup of a
superprocess with underlying spatial motion $\xi$ and branching mechanism
$(x,f)\mapsto\phi_*(f(x)) + \iittPsi_*(x,f)$. Since clearly
$\phi_*'(\lambda)\to\infty$ as $\lambda\to\infty$, the result
follows by
Theorem 8.6 of \citet{Li11}.
\end{pf}

Let us consider a canonical c\`adl\`ag realization $Y = (\iittOmega,
\mcr{G}, \mcr{G}_t, Y_t, \mbf{Q}_\mu)$ of the nonlocal branching
superprocess with transition semigroup $(Q_t)_{t\ge0}$ defined by
(\ref{6.13}) and (\ref{7.1}), where $\iittOmega= D([0,\infty),
M(T))$. Let
$Y_t(q) = Y_t[0,q]$ for $t\ge0$ and $q\in T$. By Theorem \ref{t6.4}, we
have
%
%e7.3 #&#
\begin{equation}
\label{7.3} \mbf{Q}_\mu\exp\biggl\{-\int_0^\infty
Y_s(q)f(s)\,ds \biggr\} = \exp\bigl\{-\mu[0,q]u_q(0,f)
\bigr\},
\end{equation}
where $s\mapsto u_q(s,f)$ is the unique compactly supported bounded
positive solution to (\ref{4.7}). By Theorems 8.22 and 8.23 of \citet{Li11},
for each $x\in T$ there is an excursion law $\mbf{N}_x$ on
$D([0,\infty),
M(T))$ of the superprocess such that $\mbf{N}_x\{Y_0\neq0\} = 0$ and
%
%e7.4 #&#
\begin{equation}
\label{7.4} \mbf{N}_x \bigl[1 - e^{-\int_0^\infty\langle Y_s,f_s\rangle
\,ds} \bigr] = - \log
\mbf{Q}_{\delta_x} \bigl[e^{-\int_0^\infty\langle
Y_s,f_s\rangle \,ds} \bigr]
\end{equation}
for any bounded positive Borel function $(s,y)\mapsto f_s(y)$ on
$[0,\infty)\times T$ with compact support. In view of (\ref{7.3}) and
(\ref{7.4}), for any $f\in B^+[0,\infty)$ with compact support, we have
%
%e7.5 #&#
\begin{equation}
\label{7.5} \mbf{N}_0 \bigl[1 - e^{-\int_0^\infty Y_s(q)f(s)\,ds} \bigr] =
u_q(0,f),
\end{equation}
where $s\mapsto u_q(s,f)$ is the unique compactly supported bounded
positive solution to (\ref{4.7}).
%
%th7.2 #&#
\begin{theorem}\label{t7.2} Under $\mbf{N}_0$ the path-valued process
$\{(Y_t(q))_{t\ge0}\dvtx  q\in T\}$ satisfies the Markov property with
transition semigroup $\{\mbf{P}_{p,q}\dvtx  q\ge p\in T\}$ such that
%
%e7.6 #&#
\begin{equation}
\label{7.6} \qquad\int_{D^+[0,\infty)} e^{-\int_0^\infty f(s)w(s) \,ds} \mbf{P}_{p,q}(
\eta,dw) = \exp\biggl\{-\int_0^\infty
u_{p,q}(s,f)\eta(s) \,ds \biggr\},
\end{equation}
where $f\in B^+[0,\infty)$ has compact support, and $u_{p,q}(s,f)$ is
defined by (\ref{4.8}).
\end{theorem}
\begin{pf}
By Theorem \ref{t6.4}, under $\mbf{Q}_{\delta_0}$ the process
$\{(Y_t(q))_{t\ge0}\dvtx  q\in T\}$ satisfies the Markov property with
transition semigroup defined by (\ref{7.6}). Suppose that
$(s,x)\mapsto
f_s(x)$ is a bounded positive Borel function on $[0,\infty)\times T$, and
$s\mapsto g_s$ is a bounded positive Borel function on $[0,\infty)$, both
with compact supports. Then we have
\begin{eqnarray*}
&&\mbf{Q}_{\delta_0} \biggl[\exp\biggl\{-\int_0^\infty
\bigl[\langle Y_s,f_s1_{[0,p]}\rangle+
Y_s(q)g_s\bigr]\,ds \biggr\} \biggr]
\\
&&\qquad = \mbf{Q}_{\delta_0} \biggl[\exp\biggl\{-\int_0^\infty
\bigl[\langle Y_s,f_s1_{[0,p]}\rangle+
Y_s(p)u_{p,q}(s,g)\bigr]\,ds \biggr\} \biggr].
\end{eqnarray*}
From this and (\ref{7.4}) it follows that
\begin{eqnarray*}
&&\mbf{N}_0 \biggl[1 - \exp\biggl\{-\int_0^\infty
\bigl[\langle Y_s,f_s1_{[0,p]}\rangle+
Y_s(q)g_s\bigr]\,ds \biggr\} \biggr]
\\
&&\qquad = \mbf{N}_0 \biggl[1 - \exp\biggl\{-\int_0^\infty
\bigl[\langle Y_s,f_s1_{[0,p]}\rangle+
Y_s(p)u_{p,q}(s,g)\bigr]\,ds \biggr\} \biggr].
\end{eqnarray*}
Then subtracting the quantity
\[
\mbf{N}_0 \biggl[1 - \exp\biggl\{-\int_0^\infty
\langle Y_s,f_s1_{[0,p]}\rangle \,ds \biggr\}
\biggr]
\]
from both sides, we get
\begin{eqnarray*}
&&\mbf{N}_0 \biggl[\exp\biggl\{-\int_0^\infty
\langle Y_s,f_s1_{[0,p]}\rangle \,ds \biggr\}
\\
&&\hspace*{3.4pt}\quad{} \times\biggl(1 - \exp\biggl\{-\int_0^\infty
Y_s(q)g_s \,ds \biggr\} \biggr) \biggr]
\\
&&\qquad{} = \mbf{N}_0 \biggl[\exp\biggl\{-\int_0^\infty
\langle Y_s,f_s1_{[0,p]}\rangle \,ds \biggr\}
\\
&&\hspace*{18.1pt}\qquad\quad{} \times\biggl(1 - \exp\biggl\{-\int_0^\infty
Y_s(p)u_{p,q}(s,g) \,ds \biggr\} \biggr) \biggr].
\end{eqnarray*}
A monotone class argument shows that
\begin{eqnarray*}
&&\mbf{N}_0 \biggl[F \biggl(1 - \exp\biggl\{-\int
_0^\infty Y_s(q) g_s \,ds
\biggr\} \biggr) \biggr]
\\
&&\qquad = \mbf{N}_0 \biggl[F \biggl(1 - \exp\biggl\{-\int
_0^\infty Y_s(p) u_{p,q}(s,g)
\,ds \biggr\} \biggr) \biggr]
\end{eqnarray*}
for any positive Borel function $F$ on $D([0,\infty), M(T))$, measurable
with respect to the $\sigma$-algebra generated by $\{Y_t[0,v]\dvtx  t\ge0,
0\le v\le p\}$. That implies the desired Markov property of the process
$\{(Y_t(q))_{t\ge0}\dvtx  q\in T\}$.
\end{pf}

A characterization of the finite-dimensional distributions of the
path-valued process $\{(Y_t(q))_{t\ge0}\dvtx  q\in T\}$ under the excursion
law $\mbf{N}_0$ can be given by combining (\ref{7.5}) and (\ref{7.6}).
Similarly, one can obtain characterizations of the finite-dimensional
distributions of the path-valued process $\{(Y_t(q))_{0\le
t\le\alpha}\dvtx\break 
q\in T\}$ for $\alpha>0$ and the total mass process
\[
\sigma(q):= \int_0^\infty Y_t(q) \,dt,\qquad
q\in T.
\]
The following result should be compared with Theorem 6.7 of
\citet{AbrDel12}.
%
%th7.3 #&#
\begin{theorem}\label{t7.4} Suppose that $\phi$ is a branching
mechanism such
that $\phi(\lambda)\to\infty$ as $\lambda\to\infty$. Let $\{\phi_q\dvtx  q\in
T\}$ be the admissible family defined by (\ref{5.7}). Let $\theta\in
T$ be
a strictly positive constant. Then for any positive random variable
$G(\theta)$, measurable with respect to the $\sigma$-algebra
generated by
$\{Y_t(v)\dvtx  t\ge0, v\in[0,\theta]\}$, we have
\[
\mbf{N}_0\bigl[G(\theta)|A=\theta\bigr] = \phi'\bigl(
\phi^{-1}\bigl(\phi(-\theta)\bigr)\bigr) \mbf{N}_0\bigl[G(
\theta)\sigma(\theta) 1_{\{\sigma(\theta)<\infty\}}\bigr].
\]
\end{theorem}
\begin{pf}
Based on Theorem \ref{t7.2} and the Markov property of the
path-valued process $\{(Y_t(q))_{t\ge0}\dvtx  q\in T\}$, this follows as in
the proofs of Theorems \ref{t5.3} and \ref{t5.4}.
\end{pf}

\section*{Acknowledgments}

I would like to acknowledge the Laboratory of Mathematics and Complex
Systems (Ministry of Education) for providing me the research
facilities. I am grateful to Dr. Leif D\"oring for helpful comments on
the presentation of the paper.

%suskaldyti doi

% imsref loaded by lrinkeviciute, 2012-10-23 11:08:00
% imsref loaded by lrinkeviciute, 2012-10-23 11:09:04

\printaddresses

\end{document}